\documentclass[a4paper,10pt]{article}

\usepackage[dvips]{graphics,graphicx,epsfig}
\usepackage{fancybox}
\usepackage{latexsym}
\usepackage{amssymb}
\usepackage[all]{xy}
\usepackage{color}

\newcommand{\keywordsname}{Keywords}
\newenvironment{keywords}%
  {\small
    \list{}{\labelwidth0pt
      \leftmargin0pt \rightmargin\leftmargin
      \listparindent\parindent \itemindent0pt
      \parsep0pt
      }%
    \item[\hskip\labelsep\bfseries\keywordsname:]}{\endlist}

\def\Nset{\mathrm{I\!N}}
\def\Rset{\mathrm{I\!R}}

\def\rien{\rule{0pt}{0pt}}
\def\text#1{\mbox{#1}}
\def\counterfact{>}

\begin{document}


\title{D\emph{eterministic} m\emph{odal} B\emph{ayesian} L\emph{ogic}: derive the Bayesian inference within the modal logic $T$}

\author{\begin{tabular}{c}
{\bf Fr\'ed\'eric Dambreville}\\
D\'el\'egation G\'en\'erale pour l'Armement, DGA/CEP/GIP\\
16 Bis, Avenue Prieur de la C\^ote d'Or\\
F 94114, France\\[3pt]
{\tt http://email.FredericDambreville.com}\\
{\tt http://www.FredericDambreville.com}
\end{tabular}}

\date{}
\maketitle

\begin{abstract}
In this paper a conditional logic is defined and studied.
This conditional logic, DmBL, is constructed as a deterministic counterpart to the Bayesian conditional.
The logic is unrestricted, so that any logical operations are allowed.
A notion of logical independence is also defined within the logic itself.
This logic is shown to be non-trivial and is not reduced to classical propositions.
A model is constructed for the logic.
Completeness results are proved.
It is shown that any unconditioned probability can be extended to the whole logic DmBL.
The Bayesian conditional is then recovered from the probabilistic DmBL.
At last, it is shown why DmBL is compliant with Lewis' triviality.
\end{abstract}

%
%
\section{Introduction}
Bayesian inference is a powerful principle for modeling and manipulating probabilistic information.
In many cases, Bayesian inference is considered as an optimal and legitimate rule for inferring such information.
\begin{itemize}
\item Bayesian filters for example, and their approximations by means of sequential Monte-Carlo, are typically regarded as optimal filters~\cite{arulampalam02tutorial,nadia, brehard},
\item Bayesian networks are particularly powerful tools for modeling uncertain information, since they merge logical and independence priors, for reducing the complexity of the laws, with the possibility of learning~\cite{pearlRus,murphy,Dambreville:cepomdp}.
\end{itemize}
Although Bayesian inference is an established principle, it is recalled~\cite{jaynes} that it has been disputed until the middle of the XXth century, in particular by the frequencist community.
What made the Bayesian inference established is chiefly a logical justification of the rule~\cite{cox,jaynes}.
Some convergence with the frequencist interpretation achieved this acceptation.
Cox derived the characteristics of Probability and of Bayesian conditional from hypothesized axioms about the probabilistic system, which themselves were reduced in terms of functional equations.
Typical axioms are:
\begin{itemize}
\item The operation which maps the probability of a proposition to the probability of its negation is idempotent,
\item The probability of $A\wedge B$ depends only of the probability of $A$ and of the probability of $B$ given that $A$ is true,
\item The probability of a proposition is independent of the way it is deduced (consistency).
\end{itemize}
It is noticed that Cox interpretation has been criticized recently for some imprecision and reconsidered~\cite{debrucq1,debrucq2,Halpern}.
\\[5pt]
In some sense, Cox justification of the Bayesian conditional is not entirely satisfactory, since it is implicit: it justifies the Bayesian conditional as the operator fulfilling some natural properties, but does not construct a full underlying logic  priorly to the probability.
The purpose of this paper is to construct an explicit logic for the Bayesian conditional as a conditional logic:
\begin{enumerate}
\item Build a (deterministic) conditional logic, \emph{priorly to any notion of probability}.
This logic will extend the classical propositional logic.
It will contain a conditional operator $(\cdot|\cdot)$, so that the conditional proposition $(\psi|\phi)$ could be built for any propositions $\phi$ and $\psi$,
\item Being given a probability $p$ over the unconditioned propositions, derive the probabilistic Bayesian conditional from an extension $\overline{p}$ of the probability $p$ over this conditional logic.
More precisely, the Bayesian conditional will be derived by $p(\psi|\phi)=\overline{p}\bigl((\psi|\phi)\bigr)$\,.
\end{enumerate}
The construction of an explicit underlying logic provides a better understanding of the Bayesian conditional, but will also make possible the comparison with other rules for manipulating probabilistic information, based on other logics~\cite{dambreville}.
\\[5pt]
It is known that the construction of such underlying logic is heavily constrained by Lewis' triviality \cite{lewis}, which has shown some critical issues related to the notion of conditional probability; refer also to \cite{hajek1,hajek2,fraassen}.
In particular, Lewis' result implies strong hypotheses about the nature of the conditionals.
In most cases, the conditionals have to be constructed outside the space of unconditioned propositions. 
This result implied the way the logic of Bayesian conditional has been investigated.
Many approaches do not distinguish the Bayesian conditional from probabilistic notions.
This is particularly the case of the theory called \emph{Bayesian Logic}~\cite{BayesianLogic}, which is an extension of probabilistic logic programming by the way of Bayesian conditioning.
Other approaches like conditional logic or algebra result in the construction of conditional operators, which finally arise as abstraction independent of any probability.
These logical constructions are approximating the Bayesian conditional or are constrained in use. 
The Bayesian inference is a rich notion.
It is also evocated subsequently, how the Bayesian inference has been applied to the definition of \emph{default reasoning systems}~\cite{bourne}.
\\[5pt]
Since Lewis' triviality is a fundamental reference in this work, it is introduced now.
By the way, different logical approaches of the Bayesian conditional are evocated, and it is shown how these approaches avoid the triviality.
%
\paragraph{Lewis' triviality.}
Let $\Omega$ be the set of all events, and $\mathcal{M}$ be the set of measurable subsets of $\Omega$.
Let $\mathit{Pr}(\mathcal{M})$ be the set of all probability measures defined on $\mathcal{M}$.
Lewis' triviality \cite{lewis} may be expressed as follows:
\\[5pt]
\emph{Let $A,B\in \mathcal{M}$ with $\emptyset\subsetneq B\subsetneq 
A\subsetneq \Omega$\,.
Then, it is impossible to build a proposition $(B|A)\in \mathcal{M}$ such that $\pi\bigl((B|A)\bigr)=\pi(B|A)\stackrel{\Delta}{=}\frac{\pi(A\cap B)}{\pi(A)}$ for any $\pi\in\mathit{Pr}(\mathcal{M})$ with $0<\pi(B)<\pi(A)<1$\,.}
\\[5pt]
Lewis' triviality thus makes impossible the construction of a (Bayesian) conditional operator within the same Boolean space.
\begin{description}
\item[Proof.]
For any propositions $C,D$\,, define $\pi_C(D)=\pi(D|C)=\frac{\pi(C\cap D)}{\pi(C)}$\,, when $\pi(C)>0$.\\
The proof of Lewis' result relies of the following calculus:
\begin{equation}
\label{Eq:DmBL:v2:Lewis:1}
\begin{array}{@{}l@{}}\displaystyle\vspace{5pt}
\pi((B|A)|C)=\pi_C((B|A))=\pi_C(B|A)=\frac{\pi_C(A\cap B)}{\pi_C(A)}
\\\displaystyle
\rien\hspace{100pt}=\frac{\frac{\pi(C\cap A\cap B)}{\pi(C)}}{\frac{\pi(A\cap C)}{\pi(C)}}= \frac{\pi(C\cap A\cap B)}{\pi(A\cap C)}=\pi(B|C\cap A)\;.
\end{array}
\end{equation}
Denoting $\sim B=\Omega\setminus B$, it is inferred then:
$$
\begin{array}{@{}l@{}}\displaystyle
\frac{\pi(B)}{\pi(A)}=\frac{\pi(A\cap B)}{\pi(A)}=\pi(B|A)=\pi((B|A))=\pi((B|A)|B)\pi(B)+\pi((B|A)|\sim B)\pi(\sim B)
\vspace{3pt}\\\displaystyle
\hspace{15pt}=\pi(B|B\cap A)\pi(B)+\pi(B|\sim B\cap A)\pi(\sim B)=1\times \pi(B)+0\times \pi(\sim B)=\pi(B)\;,
\end{array}
$$
which contradicts the hypotheses $\pi(B)\ne0$ and $\pi(A)\ne1$\,.
\item[$\Box\Box\Box$]\rien
\end{description}
In fact, the derivation~(\ref{Eq:DmBL:v2:Lewis:1}) relies on the hypothesis that $(B|A)\in\mathcal{M}$, which implies by definition of $(B|A)$ the  relation $\pi_C\bigl((B|A)\bigr)\pi_C(A)=\pi_C(B\cap A)$\,.
\\[5pt]
If the proposition $(B|A)$ is outside $\mathcal{M}$, it becomes necessary to build for any probability $\pi$ its extension $\overline{\pi}$ over the outside propositions; in particular, it will be defined $\overline{\pi}\bigl((B|A)\bigr)=\pi(B\cap A)/\pi(A)$, for any $A,B\in\mathcal{M}$\,.
In practice, there is no reason to have $\overline{\pi_C}(D)=\overline{\pi}\bigl((D|C)\bigr)$ for $D\not\in\mathcal{M}$\,; then, the above triviality does not work anymore.
\\[5pt]
The property $\overline{\pi_C}\ne\overline{\pi}\bigl((\cdot|C)\bigr)$ is somewhat counter-intuitive.
In particular, it means that conditionals are not conserved by conditional probabilities.
However, it allows the construction of a conditional logic for the Bayesian conditional; our work provides an example of such construction.
\paragraph{Probabilistic logic and Bayesian logic.}
Probabilistic logic, as defined by Nilsson~\cite{nilss,paass,pearl}, has been widely studied in order to model and manipulate the uncertain information.
It tracks back from the seminal work of Boole~\cite{boole}.
In probabilistic logic, the knowledge, while logically encoded by classical propositions, is expressed by means of constraints on the probability over these propositions.
For example, the knowledge over the propositions $A,B$ may be described by:
\begin{equation}
\label{Eq:DmBL:v2:BL:1}
v_1\le p(A)\le v_2\quad\mbox{and}\quad v_3\le p(A\rightarrow B)\le v_4 \;,
\end{equation}
where $A\rightarrow B\equiv\neg A\vee B$ and $v_i|1\le i\le 4$ are known bound over the probabilities.
Equations like (\ref{Eq:DmBL:v2:BL:1}) turn out to be a linear set of constraints over $p$.
In other words, it is then possible to characterize all the possible values for $p$ by means of a linear system.
For example, if our purpose is to know the possible values for $p(B)$, we just have to solve:
$$
\left\{\begin{array}{l@{}}
v_1\le p(A\wedge B)+p(A\wedge \neg B)\le v_2\;,\\
v_3\le p(A\wedge B)+p(\neg A\wedge B)+p(\neg A\wedge \neg B)\le v_4 \;,\\
p(A\wedge B)+p(\neg A\wedge B)=p(B)\;.
\end{array}\right.
$$
Notice that probabilistic logic by itself does not manipulate conditional probabilities or any notion of independence.
Proposals for extending the probabilistic logic to conditionals has appeared rather early~\cite{adams},
but Andersen and Hooker~\cite{andersen,BayesianLogic} introduced an efficient modeling and solve of such problems.
This new paradigm for manipulating Bayesian probabilistic constraints has been called \emph{Bayesian Logic}.
\\[5pt]
For example, let us introduce a new proposition $C$ to problem~(\ref{Eq:DmBL:v2:BL:1}).
Assume now that the new system is characterized by a bound over the conditional probability $p(C|A\wedge B)$ and by an independence hypothesis between $A$ and $B$.
The set of constraints could be rewritten as:
$$\left\{\begin{array}{l@{}}
v_1\le p(A)\le v_2\quad\mbox{and}\quad v_3\le p(A\rightarrow B)\le v_4 \;,
\\
v_5\le p(C|A\wedge B)\le v_6\quad\mbox{and}\quad p(A\wedge B)=p(A)p(B)\;.
\end{array}\right.$$
The constraint on $p(C|A\wedge B)$ turns out to be linear, since it could be rewritten:
$$
v_5\,p(A\wedge B)\le p(C\wedge A\wedge B)\le v_6\,p(A\wedge B)\;,
$$
but constraint $p(A\wedge B)=p(A)p(B)$ remains essentially a non-linear constraint.
Constraints involving both conditional and non-conditional probabilities also generate non-linearity.
In~\cite{andersen} and more thoroughly in~\cite{BayesianLogic}, Andersen and Hooker expose a methodology for solving these non-linear programs.
In particular, the structure of the Bayesian Network is being used in order to reduce the number of non-linear constraints.
\\[5pt]
\emph{Bayesian Logic} is a paradigm for solving probabilistic constraint programs, which involve Bayesian constraints.
Since it does not construct the Bayesian conditional as a strict logical operator, this theory is not concerned by Lewis' triviality.
\emph{Bayesian Logic} departs fundamentally from our approach, since \emph{Deterministic modal Bayesian Logic} intends to build the logic underlying the Bayesian conditional priorly to the notion of probability.
\paragraph{Conditional Event algebra.}
In Conditional Event Algebra~\cite{goodman,goodman2, dubois,calabrese}, the conditional could be seen as an external operator $(\,|\,)$\,, which maps pairs of unconditioned propositions toward an \emph{external} Boolean space, while satisfying the properties related to the Bayesian conditional.
There are numerous possible constructions of a CEA.
In fact, most CEAs provide conditional rules which are richer than the strict Bayesian conditional.
For example, the CEA, denoted DGNW~\cite{goodman}, is characterized by the following properties:
\begin{equation}
\label{DmBLv2:DGNW:eq1:0}
\left\{\begin{array}{@{}l@{}}
\vspace{5pt}\displaystyle
(a|b)\wedge(c|d)=\bigl(a\wedge b\wedge c\wedge d \big| (\neg a\wedge b)\vee(\neg c\wedge d)\vee(b\wedge d)\bigr)\;,
\\\displaystyle
(a|b)\vee(c|d)=\bigl((a\wedge b)\vee (c\wedge d) \big| (a\wedge b)\vee(c\wedge d)\vee(b\wedge d)\bigr)\;.
\end{array}\right.
\end{equation}
Property~(\ref{DmBLv2:DGNW:eq1:0}) infers the general Bayesian rule:
$$
(a|b\wedge c)\wedge (b|c)=(a\wedge b|c)\,,
$$
but also a Boolean morphism:
\begin{equation}
\label{DmBLv2:DGNW:eq1:1}
(a\wedge b|c)=(a|c)\wedge(b|c)\quad\mbox{and}\quad(a\vee b|c)=(a|c)\vee(b|c)\,.
\end{equation}
Notice that the external space hypothesis is fundamental here, and it is not possible to write $(a|\Omega)=a$ where $\Omega$ is the set of all events.
In particular, DGNW allows $(a|b)\wedge (b|\Omega)=(a\wedge b|\Omega)$, but not $(a|b)\wedge b=a\wedge b$\,.
\\[5pt]
Now, property~(\ref{DmBLv2:DGNW:eq1:0}) defines much more combinations than the strict Bayesian conditional.
Indeed, the combination $(a|b)\wedge(c|d)$ is reduced for any choice of $a,b,c,d$, which is not possible with a classical Bayesian combination.
\\[5pt]
The counterpart of such nice properties is the necessity to restrict the conditional to unconditioned propositions.
It is yet proposed in~\cite{goodman} a closure of DGNW:
\begin{equation}
\label{DmBLv2:DGNW:eq2:0}
\bigl((a|b)\big|(c|d)\bigr)=\big(a\big|(b\wedge\neg a\wedge\neg d)\vee(b\wedge c\wedge d)\bigr)\;,
\end{equation}
but this closure is not compatible with a probabilistic interpretation and \emph{fails} to satisfy the intuitive relation $P\bigl((a|b)\big|(c|d)\bigr)=P\bigl((a|b)\wedge(c|d)\bigr)/P(c|d)$\,.
More generally, CEAs are practically restricted to only one level of conditioning, and usually avoid any interferences between unconditioned and conditioned propositions.
These restrictions are also the way, by which CEAs avoid Lewis' triviality.
\paragraph{Conditional logics.}
\emph{Conditional} is an ambiguous word, since there may be different meaning owing to the community.
Even the classical inference, $\phi\rightarrow\psi\equiv \neg\phi\vee\psi$\,, is called \emph{material conditional}.
Despite classical inference is systematically used by mathematicians, its disjunctive definition makes it improper for some conditions of use.
For example, it is known that it is by essence non-constructive, an issue which tracks back to the foundation of modern mathematic~\cite{intuitionism}.
Somehow, the possibility to infer from the contradiction is also counter-intuitive; by the way, the Bayesian inference makes no sense, while inferring from the contradiction.
\\[5pt]
In the framework of Bayesian inference, we are interested in defining the law of a conditional independently of the factual state of the hypothesis.
Thus, the Bayesian conditional could be related to the notion of \emph{Counterfactual conditional}.
The classical inference is actually not counterfactual, since it is by definition dependent of the hypothesis.
From now on, the notion of conditional will refer to counterfactual conditionals, and related extensions.
\\[5pt]
A non-classical inference is often related to a modal paradigm.
While first defining counterfactual conditionals
(an example of such conditional, VCU, is detailed in section~\ref{Section:Theorem:DmBL:sub:2}), the philosophers David Lewis and Robert Stalnaker~\cite{Lewis:counterfactuals,stalnaker} based their model constructions on the \emph{possible world semantics} of modal logic (other authors also consider Kripke model extensions~\cite{giordano1}).
Stalnaker claimed that it was possible to construct such conditional, denoted $>$, within the universe of events, so as to match the probabilistic Bayesian conditional, \emph{i.e.} $p(A > B)=p(B|A)=p(A\wedge B)/p(A)$.
Lewis answered negatively~\cite{lewis} to this conjecture.
However, this was not the end of the interaction between conditional logics and probabilities.
On the basis of the semantic, Lewis proposed an alternative interpretation of the probability $p(A > B)$, called Imaging~\cite{lepage}.
\\[5pt]
It is not our purpose to detail these notions here, but we point out the following:
\begin{itemize}
\item The probabilistic interpretation of the conditional by imaging does not provide an interpretation of the probabilistic Bayesian conditional.
Nevertheless, Stalnaker's conjecture could be weakened so as to overcome the triviality:
as already explained in our previous discussion, the triviality could be avoided by constructing the conditionals outside the classical propositions space and extending the probability accordingly,
\item
The existing conditional logics are still rough approximations of the Bayesian conditional.
This point will be discussed in details later.
But for example, it seems that the negation operators of the existing conditional logics are usually relaxed, when compared to the Bayesian conditioning (refer to the logic VCU defined in section~\ref{Section:Theorem:DmBL:sub:2}).
Typically, a relation like \mbox{$\neg(\phi>\psi)\equiv\phi>\neg\psi$}\,, \emph{i.e.} a logical counterpart of $p(\neg\psi|\phi)=1-p(\psi|\phi)$\,, is not retrieved.\footnote{This property could also be related to the Boolean morphism~(\ref{DmBLv2:DGNW:eq1:1}) of CEA.}
It contradicts the axiom $\vdash\phi>\phi\mbox{ (Id)}$ which is widely accepted in the literature; refer to deduction (\ref{DmBL:VCU:eq1}) in section~\ref{Section:Theorem:DmBL:sub:2}.
\end{itemize}
\paragraph{Default reasoning.}
Default reasonings are related to logical systems which are able to deduce what \emph{normally} happens, when only partial information are available.
The idea is to use default rules or default information.
Then, the most plausible assumptions are deduced, in regards to the current information.
Of course, some conclusions may be retracted, if they are corrected by new sources of information.
\\[5pt]
It is noticed that the Bayesian inference is able of some kind of default reasoning, by adapting the belief of a proposition according to an hypothesis.
In fact, the system $P$ for default reasoning could be derived from a Bayesian interpretation~\cite{adams,definetti,bourne}.
The $\varepsilon-$semantic is particularly enlightening.
Let us denote $\phi\Rightarrow\psi$ the default, which means ``$\psi$ is normally true, when $\phi$ holds true''.
Then let us interpret $\phi\Rightarrow\psi$ by the constraint $p(\psi|\phi)\ge1-\varepsilon$, where $\varepsilon$ is an infinitesimal.\footnote{Notice that the definition of $\Rightarrow$ is a meta-definition (or also a second-order definition)\,.}
Now, it is deduced from the Bayesian inference:
\begin{equation}
\label{DmBL:Def:Eq:1}
p(\eta|\phi\wedge\psi)=\frac{p(\eta|\phi)-\bigr(1-p(\psi|\phi)\bigl)p(\eta|\phi\wedge\neg\psi)}{p(\psi|\phi)}\,.
\end{equation}
From hypotheses $\phi\Rightarrow\psi$ and $\phi\Rightarrow\eta$, it would come $p(\psi|\phi)\ge1-\varepsilon_1$ and $p(\eta|\phi)\ge1-\varepsilon_2$, where $\varepsilon_1$ and $\varepsilon_2$
are infinitesimals, and then by~(\ref{DmBL:Def:Eq:1}):
$$
p(\eta|\phi\wedge\psi)\ge1-\epsilon_1-\epsilon_2
\quad\mbox{and}\quad
\phi\wedge\psi\Rightarrow\eta\;.
$$
This deduction conduces to the well known \emph{Cautious Monotonicity} rule:
$$
CM:\ \frac{\phi\Rightarrow\psi,\ \phi\Rightarrow\eta}{\phi\wedge\psi\Rightarrow\eta}\;.
$$
The rule CM already infers some default reasoning:
\begin{quote}
Assume $b\Rightarrow f$ (birds normally fly), $p\Rightarrow b$ (penguins normally are birds) and $p\Rightarrow \neg f$ (penguins normally do not fly).
Then by applying CM on $p\Rightarrow b$ and $p\Rightarrow \neg f$, it comes $p\wedge b\Rightarrow \neg f$\,.
The \emph{non-monotonic} inference $\Rightarrow$ is thus able to handle sub-cases.
\end{quote}
From the infinitesimal probabilistic interpretation, the following rules of $P$ are also deduced:
$$\begin{array}{lllll}
~\mbox{\footnotesize Id:}~\frac{}{\phi\Rightarrow\phi}
&
~\mbox{\footnotesize RW:}~\frac{\vdash\psi\rightarrow\eta,\ \phi\Rightarrow\psi}{\phi\Rightarrow\eta}
&
\frac{\vdash\phi\leftrightarrow\psi,\ \phi\Rightarrow\eta}{\psi\Rightarrow\eta}
&
\frac{\phi\Rightarrow\psi,\ \phi\Rightarrow\eta}{\phi\Rightarrow\psi\wedge\eta}
&
\frac{\phi\Rightarrow\eta,\ \psi\Rightarrow\eta}{\phi\vee\psi\Rightarrow\eta}
\end{array}\;.$$
System $P$ shares some common rules with existing conditional logics (\emph{e.g.} (Id) is in VCU and \emph{Right Weakening} (RW) is a sub-case of rule (CR) of VCU; refer to section~\ref{Section:Theorem:DmBL:sub:2}).
However, the default $\Rightarrow$ appears as a meta-operator, which acts as an external operator over the classical propositions only.
Contrary to CEA DGNW, equation~(\ref{DmBLv2:DGNW:eq1:1}), the default propositions $\phi\Rightarrow\psi$ do not constitute a Boolean space which maps from the classical propositions by a Boolean morphism.
Both properties (internal operator and Boolean morphism) are desirable, for a logical interpretation of the Bayesian conditional.
But system $P$ allows default reasoning, while DGNW and VCU (and DmBL, defined subsequently) are incompatible with it (equation~(\ref{DmBLv2:DGNW:eq1:0}), axiom Ax.4 of VCU defined in section~\ref{Section:Theorem:DmBL:sub:2}).
\\[5pt]
Also related to the infinitesimal interpretation is the ranking of the models (typically, a rank is interpreted as an infinitesimal order).
The ranking conditions the way the defaults are prioritized.
There has been various extension of system $P$ (alternative inferences, ranking methods).
In some recent works, Lukasiewicz et al.~\cite{biazzo} proposed a probabilistic extension of default reasonings (including system $P$), which is not restricted to an infinitesimal probabilistic interpretation.
Each default is then associated to an interval of constraint for its probability.
\paragraph{Why a new conditional logic?}
The previous approaches and uses of the Bayesian logic imply restrictions or approximations to the logical interpretation of the Bayesian conditional.
\emph{Bayesian logic} does not provide a logical interpretation of the Bayesian conditional, but rather a methodology for solving the program related to a probabilistic Bayesian modeling.
\emph{Conditional event algebras} provide an interesting logical interpretation of the Bayesian conditional, but are highly constrained in their definition;
from a logical viewpoint, the impossibility to handle or to combine multiple-levels conditionals
constitutes a limitation in terms of the coherence of the models.
Existing \emph{conditional logics} are insufficient for characterizing the Bayesian conditional properly.
The default reasoning systems and their interpretation by means of the probabilistic Bayesian inference are quite interesting and a source of inspiration, but they are mainly a consequence of the  probabilistic Bayesian inference and account for the Bayesian conditional only partially.
Our work intends to overcome these limitations, by constructing a new conditional logic which is in accordance with the Bayesian conditional.
It is not our purpose, however, to enrich the Bayesian rule, as it is proposed in most CEA.
\\[5pt]
Our logic, denoted \emph{Deterministic modal Bayesian Logic} (DmBL), is constructed according to a modal background (system T).
The conditional operator is defined in parallel to a relation of logical independence. This relation is defined within the logic, and not at a meta-level.
The probabilistic Bayesian inference is recovered from the derived logical theorems and the logical independence.
This process implies an extension of probability from the unconditioned logic toward DmBL.
As a final result, a theorem is proved that guarantees the existence of such extension (Lewis' result is thus avoided).
\\[5pt]
Section~\ref{Section:Def:DmBL} is dedicated to the definition of the Deterministic modal Bayesian Logic.
The languages, axioms and rules are introduced.
In section~\ref{Section:Theorem:DmBL}, several theorems of the logic are derived.
A purely logical interpretation of Lewis' triviality is made, and DmBL is compared with other systems. 
A model for DmBL is constructed in section~\ref{Section:Model:DmBL}.
A partial completeness theorem is derived.
The extension of probabilities over DmBL is investigated in section~\ref{Section:Proba:DmBL}.
The probabilistic Bayesian inference is recovered from this extension.
The paper is then concluded.
%
%
%
%
%
\section{Definition of the logics}
\label{Section:Def:DmBL}
D\emph{eterministic} m\emph{odal} B\emph{ayesian} L\emph{ogic} was first defined without modality as D\emph{eterministic} B\emph{ayesian} L\emph{ogic} in a previous version of this document~\cite{Dambreville:DBL}.
The non-modal definition is uneasy to handle.
For this new modal definition, we have been inspired by the seminal work of Lewis \cite{Lewis:counterfactuals}, and also by more recent works of Laura Giordano et al., which use the modality for specifying the conditional behavior~\cite{giordano1,giordano2}.
The use of the modality is instrumental here; for this reason, m\emph{odal} appears in lower case in our terminology.
\\[5pt]
Modal logic is a powerful tool, but not intuitive at first sight.
We thus decided to introduce the modality softly, by interpreting it (in part) from probabilistic considerations: modalities will be used in order to characterize properties generally true for any possible probability.
Nevertheless, this instrumental use of modalities allows an abstraction which makes our construction independent to any notion of probability.
%
\paragraph{Introducing the modal notation.}
This paragraph intends to explain the intuition behind the subsequent modal definition of DmBL.
Some generalizations are not justified here,
but the axioms are extrapolated from results, which are easily proved.
\\[5pt]
The logic of a system is the collection of behaviors which are common to any instance of this system.
Let us consider the example of probability on a finite \emph{(unconditioned)} propositional space.
For convenience, define $\mathbb{P}$ the set of strictly positive probabilities over this space, that is $p\in\mathbb{P}$ is such that $p(\phi)>0$ for any non-empty proposition~$\phi$\,:
$$
\mathbb{P}=\bigl\{p\;\big/\;p\mbox{ is a probability and }\forall\phi\not\equiv\bot,\,p(\phi)>0\bigr\}\,.
$$ 
When the proposition $\phi$ is always \emph{true} (\emph{i.e.} $\phi$ is the \emph{set of all events}), it is known that $p(\phi)=1$ for any possible probability $p$.
This could be interpreted logically as follows:
\begin{equation}
\label{def:intro:1}
\vdash\phi\quad\mbox{implies}\quad\vdash\forall p\in\mathbb{P},\;p(\phi)=1\,.
\end{equation}
This is a typical logical relation related to the probabilities.
Notice however that:
$$
\phi\rightarrow\forall p\in\mathbb{P},\;p(\phi)=1\quad\mbox{is false in general}\,,
$$
because unless $\phi$ is always true, the property $\forall p\in\mathbb{P},\;p(\phi)=1$ is necessary false.\footnote{In first order logic, being able to prove $B$ when hypothesizing $A$ implies that $A\rightarrow B$ is proved (assumption discharge).
On this second order example, the assumption discharge does not work anymore!}
\\
It is also obvious that $\forall p\in\mathbb{P},\;p(\phi)=1$ infers $\phi$ (consider the cases $\vdash\phi$ and $\nvdash\phi$)\,:
\begin{equation}
\label{def:intro:2}
\vdash(\forall p\in\mathbb{P},\;p(\phi)=1)\rightarrow\phi\,.
\end{equation}
Another proposition is easily derived:
\begin{equation}
\label{def:intro:3}
\vdash(\forall p\in\mathbb{P},\;p(\neg\phi\vee\psi)=1)\rightarrow\bigl((\forall p\in\mathbb{P},\;p(\phi)=1)\rightarrow (\forall p\in\mathbb{P},\;p(\psi)=1)\bigr)\,.
\end{equation}
This last proposition could be considered as a \emph{modus ponens} encoded by means of the probabilities (it is recalled that $\phi\rightarrow\psi\equiv\neg\phi\vee\psi$).
Now, by using the abbreviation $\Box\phi=\bigl(\forall p\in\mathbb{P},\;p(\phi)=1\bigr)$, the propositions (\ref{def:intro:1}),  (\ref{def:intro:2}) and (\ref{def:intro:3}) are turned into:
$$\begin{array}{@{}l@{}}
\vdash\phi\quad\mbox{implies}\quad\vdash\Box\phi\,,
\\
\vdash\Box\phi\rightarrow\phi\,,
\\
\vdash\Box(\phi\rightarrow\psi)\rightarrow(\Box\phi\rightarrow \Box\psi)\,.
\end{array}$$
These are exactly the modal axioms and rule of the system~T of modal logic.
System~T is the backbone of our conditional logic, DmBL.
Additional axioms are also introduced for characterizing the conditional $(\cdot|\cdot)$.
These axioms are extrapolated from unconditioned probabilistic characterizations.
For example, it is proved for unconditioned propositions:
\begin{equation}
\label{def:intro:4:-1:-1}
\forall p\in\mathbb{P},\,p(\phi)+p(\psi)=1\quad\mbox{implies}\quad\phi\equiv\neg\psi\;.
\end{equation}
Since it is also known that $p(\psi|\phi)+p(\neg\psi|\phi)=1$ for any $p\in\mathbb{P}$, it will be assumed the axiom:
\begin{equation}
\label{def:intro:4:-1:-1:1}
(\psi|\phi)\equiv\neg(\neg\psi|\phi)\;.
\end{equation}
Of course, property~(\ref{def:intro:4:-1:-1}) normally holds for unconditioned propositions only, so that axiom~(\ref{def:intro:4:-1:-1:1}) comes in fact from an extrapolation of~(\ref{def:intro:4:-1:-1}).
\\[5pt]
Similarly, it is noticed that:
\begin{equation}
\label{def:intro:4:-1}
\forall p\in\mathbb{P},\,p(\phi)+p(\psi)\le p(\eta)+p(\zeta)\quad\mbox{implies}\quad\vdash(\phi\vee\psi)\rightarrow(\eta\vee\zeta)\;.
\end{equation}
Since $p(\psi|\phi)+p(\bot)\le p(\neg\phi\vee\psi)+p(\bot)$ for any $p\in\mathbb{P}$, it will be assumed the axiom:
\begin{equation}
\label{def:intro:4:-1:1}
\vdash(\psi|\phi)\rightarrow(\phi\rightarrow\psi)\;.
\end{equation}
Since $p(\psi\vee\eta|\phi)+p(\bot)\le p(\psi|\phi)+p(\eta|\phi)$ for any $p\in\mathbb{P}$, it is also extrapolated that
\mbox{$\vdash(\psi\vee\eta|\phi)\rightarrow\bigl((\psi|\phi)\vee(\eta|\phi)\bigr)\,.$}
Then, by applying~(\ref{def:intro:4:-1:-1:1}), it comes:
\begin{equation}
\label{def:intro:4:-1:2}
\vdash(\psi\rightarrow\eta|\phi)\rightarrow\bigl((\psi|\phi)\rightarrow(\eta|\phi)\bigr)\;,
\end{equation}
which constitutes a modus ponens for the conditional.
Axioms (\ref{def:intro:4:-1:-1:1}), (\ref{def:intro:4:-1:1}) and (\ref{def:intro:4:-1:2}) are not completely new.
In particular, they infer the Boolean morphism~(\ref{DmBLv2:DGNW:eq1:1}) of CEA.
They are not fully implemented by the existing conditional logics, however.
\\[5pt]
In order to complete this introduction, it is also noticed that:
\begin{equation}
\label{def:intro:4}
\vdash\bigl(\forall p\in\mathbb{P},\;p(\neg\phi\vee\psi)=1\bigr)\rightarrow \bigl((\forall p\in\mathbb{P},\;p(\phi)=0) \vee (\forall p\in\mathbb{P},\;p(\psi|\phi)=1)\bigr)\,.
\end{equation}
The interpretation and proof of~(\ref{def:intro:4}) is simple: when $\phi$ is a ``subset'' of $\psi$, then either $\phi$ is empty or $p(\psi|\phi)=1$ for any strictly positive probability $p$.
From~(\ref{def:intro:4}), it is then extrapolated:
\begin{equation}
\label{def:intro:5}
\vdash\Box(\phi\rightarrow\psi)\rightarrow \bigl(\Box\neg\phi \vee \Box(\psi|\phi)\bigr)\,.
\end{equation}
In fact, these extrapolated axioms imply constraints, when extending the probabilities $p\in\mathbb{P}$ over the conditioned propositions.
This paper intends to prove that these constraints are actually valid, in regard to Lewis' triviality.
It is now time for the logic definition.
\paragraph{Language.}
Let $\Theta=\{\theta_i/i\in I\}$ be a set of atomic propositions.
\\[5pt]
The language $\mathcal{L}_C$ of the classical logic related to $\Theta$ is the smallest set such that:
$$\left\{\begin{array}{@{\,}l@{}}
\Theta\subset\mathcal{L}_C
\\[5pt]
\neg\phi\in\mathcal{L}_C \mbox{ and }\phi\rightarrow\psi\in\mathcal{L}_C
\mbox{ for any }\phi,\psi\in\mathcal{L}_C
\end{array}\right.$$
The language $\mathcal{L}_T$ of the modal logic $T$ related to $\Theta$ is the smallest set such that:
$$\left\{\begin{array}{@{\,}l@{}}
\Theta\subset\mathcal{L}_T
\\[5pt]
\neg\phi\in\mathcal{L}_T\;,\ \Box\phi\in\mathcal{L}_T \mbox{ and }\phi\rightarrow\psi\in\mathcal{L}_T
\mbox{ for any }\phi,\psi\in\mathcal{L}_T
\end{array}\right.$$
The language $\mathcal{L}$ of the D\emph{eterministic} m\emph{odal} B\emph{ayesian} L\emph{ogic} related to $\Theta$ is the smallest set such that:
$$\left\{\begin{array}{@{\,}l@{}}
\Theta\subset\mathcal{L}
\\[5pt]
\neg\phi\in\mathcal{L}\;,\ \Box\phi\in\mathcal{L}\;,\ \phi\rightarrow\psi\in\mathcal{L}\mbox{ and }(\psi|\phi)\in\mathcal{L}
\mbox{ for any }\phi,\psi\in\mathcal{L}
\end{array}\right.$$
In the construction of the propositions, the unary operators have priority over the binary operators; for example $\neg\phi\vee\psi=(\neg\phi)\vee\psi$\,.
The following abbreviations are defined:
\begin{itemize}
\item
$\phi\vee\psi=\neg\phi\rightarrow\psi$\,,\quad $\phi\wedge\psi=\neg(\neg\phi\vee\neg\psi)$\quad
and\quad
$\phi\leftrightarrow\psi=(\phi\rightarrow\psi)\wedge(\psi\rightarrow\phi)$\,,
\item It is chosen a proposition $\theta\in\Theta$, and it is then denoted $\top=\theta\rightarrow\theta$ and $\bot=\neg\top$\,,
\item $\Diamond\phi=\neg\Box\neg\phi$\,,
\item $\psi\times\phi=\Box\bigl((\psi|\phi)\leftrightarrow\psi\bigr)$\,.
\end{itemize}
$\top$ and $\bot$ are idealistic notations for the tautology and the contradiction.
The operator $\times$ describes the \emph{logical} independence between propositions.
The independence relation $\times$ and the conditional $(|)$ are thus conjointly defined.
\\[5pt]
We also define abbreviations for the notion of proof:
\begin{itemize}
\item $\vdash\phi$ means ``$\phi$ is proved''\,,
\item $\phi\equiv\psi$ means $\vdash\phi\leftrightarrow\psi$.
\end{itemize}
The meta-relation $\equiv$ is the logical equivalence.
\paragraph{Rules and axioms.}
The \emph{classical Logic $C$} is characterized by the \emph{Modus ponens} and the classical axioms $c\ast$ described subsequently.
\\
The \emph{modal Logic $T$} is characterized by the \emph{Modus ponens}, the classical axioms $c\ast$ and the modal rule/axioms $m\ast$ described subsequently (\emph{c.f.} also \cite{blackburn}).
\\
The D\emph{eterministic} m\emph{odal} B\emph{ayesian} L\emph{ogic}\,, \emph{i.e.} DmBL, is characterized by the \emph{Modus ponens}, the classical axioms $c\ast$, the modal rule/axioms $m\ast$ and the Bayesian axioms $b\ast$\,:
\begin{description}
\item[\rien$\quad$c1.] $\vdash \phi\rightarrow(\psi\rightarrow\phi)$\,,
\item[\rien$\quad$c2.] $\vdash (\eta\rightarrow(\phi\rightarrow\psi))\rightarrow((\eta\rightarrow\phi)\rightarrow(\eta\rightarrow\psi))$\,,
\item[\rien$\quad$c3.] $\vdash (\neg\phi\rightarrow\neg\psi)\rightarrow((\neg\phi\rightarrow\psi)\rightarrow\phi)$\,,
\item[\rien$\quad$Modus ponens.]$\vdash\phi$ and $\vdash\phi\rightarrow\psi$ implies $\vdash\psi$\,,
\item[\rien$\quad$m1.] $\vdash\phi$ implies $\vdash \Box\phi$\,,
\item[\rien$\quad$m2.] $\vdash\Box(\phi\rightarrow\psi)\rightarrow(\Box\phi\rightarrow\Box\psi)$\,,
\item[\rien$\quad$m3.] $\vdash\Box\phi\rightarrow\phi$\,,
\item[\rien$\quad$b1.] $\vdash\Box(\phi\rightarrow\psi)\rightarrow \bigl(\Box\neg\phi \vee \Box(\psi|\phi)\bigr)$\,,
\item[\rien$\quad$b2.] $\vdash(\psi\rightarrow\eta|\phi)\rightarrow\bigl((\psi|\phi)\rightarrow(\eta|\phi)\bigr)$\,,
\item[\rien$\quad$b3.] $\vdash(\psi|\phi)\rightarrow(\phi\rightarrow\psi)$\,,
\item[\rien$\quad$b4.] $\neg(\neg\psi|\phi)\equiv(\psi|\phi)$\,,
\item[\rien$\quad$b5.] \emph{($\times$ is symmetric)}~:
$\psi\times\phi\equiv\phi\times\psi$\,,
\end{description}
%
DmBL$_\ast$, a weakened version of DmBL, is defined by replacing $b5$ by the alternative axioms:
$$\begin{array}{@{}l@{}}
\mbox{\bf b5.weak.A. } \psi\times\neg\phi\equiv\psi\times\phi\;,
\\[3pt]
\mbox{\bf b5.weak.B. } \vdash\Box(\psi\leftrightarrow\eta)\rightarrow\Box\bigl((\phi|\psi)\leftrightarrow(\phi|\eta)\bigr)\;.
\end{array}$$
The axioms $m\ast$ and $b1$ to $b4$ have been introduced in the previous paragraph.
The axiom $b5$ implements the symmetry of the logical independence.
The specific notations $\vdash_C$, $\vdash_T$ and $\vdash$ will be used for denoting a proof in $C$, $T$ or DmBL/DmBL$_\ast$ respectively.
%
The following section studies the logical consequences of the axioms.
\section{Logical theorems and comparison with other systems}
\label{Section:Theorem:DmBL}
DmBL/DmBL$_\ast$ implies the classical and the T-system tautologies;
the properties of classical logic and of the T-system are assumed without proof.
Since both DmBL and DmBL$_\ast$ are studied, the possibly needed axioms $b5\ast$ are indicated in bracket.
\subsection{Theorems}
\label{Section:Theorem:DmBL:sub:1}
\emph{The proofs are done in appendix~\ref{proof:logth}.}
Next theorem is proved here as an example. 
%
\paragraph{The full universe.} $\vdash\Box\phi\rightarrow(\psi\times\phi)$\,.
In particular $(\psi|\top)\equiv\psi$\,.\\[5pt]
Interpretation: a tautology is independent with any other proposition and its sub-universe is the whole universe.
\begin{description}
\item[Proof.]
From axiom b3, it comes $\vdash(\psi|\phi)\rightarrow(\phi\rightarrow\psi)$ and  $\vdash(\neg\psi|\phi)\rightarrow(\phi\rightarrow\neg\psi)$\,.\\
Then $\vdash\phi\rightarrow\bigl((\psi|\phi)\rightarrow\psi\bigr)$ and  $\vdash\phi\rightarrow\bigl((\neg\psi|\phi)\rightarrow\neg\psi\bigr)$\,.\\
Applying b4 yields $\vdash\phi\rightarrow\bigl((\psi|\phi)\leftrightarrow\psi\bigr)$.\\
It follows $\vdash\Box\phi\rightarrow\Box\bigl((\psi|\phi)\leftrightarrow\psi\bigr)$\,.
\\[5pt]
The remaining proof is obvious.
\item[$\Box\Box\Box$]\rien
\end{description}
\paragraph{Axioms order.} Axiom b5 implies b5.weak.A.
\paragraph{The empty universe [b5.weak.A].} $\vdash\Box\neg\phi\rightarrow(\psi\times\phi)$\,.
In particular $(\psi|\bot)\equiv\psi$\,.
\paragraph{Left equivalences.}
$\vdash\Box(\psi\leftrightarrow\eta)\rightarrow\Bigl(\Box\neg\phi\vee\Box\bigl((\psi|\phi)\leftrightarrow(\eta|\phi)\bigr)\Bigr)$\,.
\\[5pt]
\emph{Corollary [b5.weak.A].} $\vdash\Box(\psi\leftrightarrow\eta)\rightarrow\Box\bigl((\psi|\phi)\leftrightarrow(\eta|\phi)\bigr)$.
\\[5pt]
\emph{Corollary 2 [b5.weak.A].} $\psi\equiv\eta$ implies $(\psi|\phi)\equiv(\eta|\phi)$.
\\
Proof is immediate from corollary.
\paragraph{Sub-universes are classical [b5.weak.A].}
\begin{itemize}
\item $(\neg\psi|\phi)\equiv\neg(\psi|\phi)$\,,
\item $(\psi\wedge\eta|\phi)\equiv(\psi|\phi)\wedge(\eta|\phi)$\,,
\item $(\psi\vee\eta|\phi)\equiv(\psi|\phi)\vee(\eta|\phi)$\,,
\item $(\psi\rightarrow\eta|\phi)\equiv(\psi|\phi)\rightarrow(\eta|\phi)$\,.
\end{itemize}
\paragraph{Evaluating $(\top|\cdot)$ and $(\bot|\cdot)$ [b5.weak.A].} Is proved $\vdash\Box\psi\rightarrow\Box(\psi|\phi)$\,.
In particular $(\top|\phi)\equiv\top$ and $(\bot|\phi)\equiv\bot$\,.
\paragraph{Inference property.}
$(\psi|\phi)\wedge\phi\equiv\phi\wedge\psi$\,.
\paragraph{Introspection.} $\vdash\Box\neg\phi\vee\Box(\phi|\phi)$\,.\\[5pt]
Interpretation: a non-empty proposition sees itself as ever true. \\
Notice that this property is compliant with $(\bot|\bot)\equiv\bot$\,.
\paragraph{Inter-independence [b5.weak.A].} $\vdash(\psi|\phi)\times\phi$\,.\\[5pt]
Interpretation: a proposition is independent of its sub-universe.
\paragraph{Independence invariance [b5.weak.A].}
$$
\begin{array}{@{}l@{}}
\vdash (\psi\times\phi)\rightarrow(\neg\psi\times\phi)\;,
\\
\vdash \bigl((\psi\times\phi)\wedge(\eta\times\phi)\bigr)\rightarrow\bigl((\psi\wedge\eta)\times\phi\bigr)\;,
\\
\vdash \Box(\psi\leftrightarrow\eta)\rightarrow\bigl((\psi\times\phi)\leftrightarrow(\eta\times\phi)\bigr)\;.
\end{array}
$$
\paragraph{Narcissistic independence.} $\vdash(\phi\times\phi)\rightarrow(\Box\neg\phi\vee\Box\phi)$\,.
\\[5pt]
Interpretation: a proposition independent with itself is either a tautology or a contradiction.
\paragraph{Independence and proof [b5.weak.A].}
$\vdash(\psi\times\phi)\rightarrow\bigl(\Box(\phi\vee\psi)\rightarrow(\Box\phi\vee\Box\psi)\bigr)$\,.\\[5pt]
Interpretation:
when propositions are independent and their disjunctions are sure, then at least one proposition is sure.
\paragraph{Independence and regularity [b5.weak.A].}
$$
\vdash\bigl((\phi\times\eta)\wedge(\psi\times\eta)\bigr)\rightarrow\Bigl(\Box\bigl((\phi\wedge\eta)\rightarrow(\psi\wedge\eta)\bigr)\rightarrow\bigl(\Box\neg\eta\vee\Box(\phi\rightarrow\psi)\bigr)\Bigr)\;.
$$
Interpretation: unless it is empty, a proposition may be removed from a logical equation, when it appears in the both sides and is independent with the equation components.
\\[5pt]
\emph{Corollary.} $\vdash\phi\times\eta$\,, $\vdash\psi\times\eta$\,, $\vdash\Diamond\eta$ and $\phi\wedge\eta\equiv\psi\wedge\eta$ implies $\phi\equiv\psi$\,.
\\[5pt]
\emph{Corollary 2.} Being given $\psi$ and $\phi$ such that $\vdash\Diamond\phi$, the proposition $(\psi|\phi)$ is uniquely defined as the solution of equation $X\wedge\phi\equiv\psi\wedge\phi$ (with unknown $X$) which is independent of $\phi$.
%
\paragraph{Right equivalences [b5].}
$\vdash\Box(\psi\leftrightarrow\eta)\rightarrow\Box\bigl((\phi|\psi)\leftrightarrow(\phi|\eta)\bigr)$ (proved with b5 but without b5.weak.B).
\\[5pt]
Interpretation: equivalence is compliant with the conditioning.
\\[5pt]
\emph{Corollary.} Axiom b5 implies b5.weak.B.
In particular, DmBL$_\ast$ is weaker than DmBL.
\\[5pt] 
\emph{Corollary of b5 or b5.weak.B.} $\psi\equiv\eta$ implies $(\phi|\psi)\equiv(\phi|\eta)$.
\\[5pt]
Together with the properties of the system $T$ and \emph{left equivalences}, this last result implies that the equivalence relation $\equiv$ is compliant with the logical operators of DmBL/DmBL$_\ast$.
In particular, replacing a sub-proposition with an equivalent sub-proposition within a theorem still makes a theorem.
\paragraph{Reduction rule [b5].}
Axiom b5 implies $\bigl(\phi\big|(\psi|\phi)\bigr)\equiv\phi$\,.
\paragraph{Markov Property [b5].}
$$
\vdash\left(\left(
\bigwedge_{\tau=1}^{t-2}\bigl((\phi_t|\phi_{t-1})\times\phi_{\tau}\bigr)
\right)\wedge\Diamond\left(\bigwedge_{\tau=1}^{t-1}\phi_{\tau}\right)\right)\longrightarrow\Box\left(
(\phi_t|\phi_{t-1})\leftrightarrow\left(\phi_t\left|\bigwedge_{\tau=1}^{t-1}\phi_{\tau}\right.\right)
\right)\;.
$$
Interpretation: the Markov property holds, when the conditioning is independent of the past and the past is possible.
\paragraph{Link between $\bigl((\eta|\psi)\big|\phi\bigr)$ and $(\eta|\phi\wedge\psi)$ [b5].}\rien\\
It is derived:
$
\bigl((\eta|\psi)\big|\phi\bigr)\wedge\phi\wedge\psi\equiv(\eta|\psi)\wedge\phi\wedge\psi\equiv\phi\wedge\psi\wedge\eta\equiv(\eta|\phi\wedge\psi)\wedge(\phi\wedge\psi)\;.
$\\[5pt]
This is a quite limited result and it is \emph{tempting} to assume the additional axiom ``$\bigl((\eta|\psi)\big|\phi\bigr)\equiv(\eta|\phi\wedge\psi)\quad\mbox{\small$(\ast)$}$''\,.
There is a really critical point here, since axiom~$(\ast)$ implies actually a logical counterpart to Lewis' triviality\,:
\begin{quote}
\emph{Let $\bigl((\eta|\psi)\big|\phi\bigr)\equiv(\eta|\phi\wedge\psi)\quad\mbox{\small$(\ast)$}$ be assumed as an axiom.\\
Then $\vdash\Diamond(\phi\wedge\psi)\rightarrow
\bigl(\Box(\phi\leftrightarrow\psi)\vee(\phi\times\psi)\bigr)$\,.}
\end{quote}
Interpretation: if $\phi$ and $\psi$ are not exclusive and not equivalent, then they are independent.
This is irrelevant and forbids the use of axiom~$(\ast)$.
\subsection{Some comparisons with other systems}
\label{Section:Theorem:DmBL:sub:2}
\paragraph{Conditional Event Algebra.}
As explained in introduction, CEAs characterize the conditional by means of an external operator.
The conditional are limited to only one level of conditioning, but consequently, the combination rules of the conditionals are richer than for the strict Bayesian conditionals.
On the contrary, DmBL handles several level of conditioning, but is not addressed to provide richer combinations than the strict Bayesian conditionals.
Thus, DmBL implements the necessary Bayesian properties:
\begin{itemize}
\item In DmBL, the conditioning constitutes a Boolean morphism (sub-universes are classical).
The property~(\ref{DmBLv2:DGNW:eq1:1}) of the CEA \emph{DGNW}~\cite{goodman} is thus recovered and generalized,
\item The Bayesian inference $p(\psi|\phi)p(\phi)=p(\phi\wedge\psi)$ is derived from the extension of probability (section \ref{Section:Proba:DmBL}) by means of the theorems $(\psi|\phi)\wedge\phi\equiv\phi\wedge\psi$ and $\vdash(\psi|\phi)\times\phi$\,.
\end{itemize}
DGNW also provides the relation \mbox{$(a|b\wedge c)\wedge (b|c)=(a\wedge b|c)\,,$}
which is related to the general Bayesian inference $p(a|b\wedge c)p(b|c)=p(a\wedge b|c)\,.$
The general Bayesian inference is a direct consequence of the Bayesian inference, and thus can be derived from DmBL too.
%
But is there a logical counterpart in DmBL to the general Bayesian inference?
This logical counterpart would be expressed by means of a double proposition:
$$
\vdash(\eta|\phi\wedge\psi)\times (\psi|\phi)
\quad\mbox{and}\quad
(\eta|\phi\wedge\psi)\wedge (\psi|\phi)\equiv (\psi\wedge\eta|\phi)\;.
$$
At this time, we are not able to decide if these propositions are derived from DmBL or are even compatible with DmBL.
\paragraph{Comparison with an existing conditional logic.}
The axioms of the conditional logic VCU (VCU is an abbreviation for the axioms system)~\cite{Lewis:counterfactuals} are considered here and compared to DmBL.
This example is representative of the difference with the other conditional logics.
\emph{Theorems derived in section~\ref{Section:Theorem:DmBL:sub:1} are referred to.}
\\[5pt]
Axioms and rules of VCU:
\\
(Ax.1) $\phi\counterfact\phi$ has a partial counterpart in DmBL, \emph{i.e.} $\vdash\Box\neg\phi\vee\Box(\phi|\phi)$ (theorem).\\
(Ax.2) $(\neg\phi\counterfact\phi)\rightarrow(\psi\counterfact\phi)$ becomes $\vdash(\phi|\neg\phi)\rightarrow(\phi|\psi)$ (derived from theorems).\\
(Ax.3) $(\phi\counterfact\neg\psi)\vee(((\phi\wedge\psi)\counterfact\xi)\leftrightarrow(\phi\counterfact(\psi\rightarrow\xi)))$ has no obvious counterpart in DmBL.\\
(Ax.4) $(\phi\counterfact\psi)\rightarrow(\phi\rightarrow\psi)$ is exactly b3.\\
(Ax.5) $(\phi\wedge\psi)\rightarrow(\phi\counterfact\psi)$ is a subcase of $\phi\wedge\psi\equiv\phi\wedge(\psi|\phi)$ (inference theorem).\\
(Ax.6) $(\neg\phi\counterfact\phi)\rightarrow\bigl(\neg(\neg\phi\counterfact\phi)\counterfact(\neg\phi\counterfact\phi)\bigr)$
becomes
$\vdash(\phi|\neg\phi)\rightarrow\bigl((\phi|\neg\phi)\big|\neg(\phi|\neg\phi)\bigr)$
(derived from theorems).\\
(CR) Counterfactual rule.
This is a multiple-task rule.
First, it allows the introduction of tautologies inside a conditional,
secondly, it implies some linearity of the conditional with $\wedge$\,:
\\\rien\hspace{50pt}Being proved $(\xi_1\wedge\dots\wedge\xi_n)\rightarrow\psi$\,, it is proved $((\phi\counterfact\xi_1)\wedge\dots\wedge(\phi\counterfact\xi_n))\rightarrow(\phi\counterfact\psi)$\,.\\
This rule is recovered in DmBL from the fact that \emph{sub-universes are classical}:
\\\rien\hspace{50pt}$\vdash(\xi_1\wedge\dots\wedge\xi_n)\rightarrow\psi$ implies $\vdash((\xi_1|\phi)\wedge\dots\wedge(\xi_n|\phi))\rightarrow(\psi|\phi)$\,.
\\[5pt]
It is noteworthy that Ax.1 and CR, with $n=1$ and $\xi_1=\phi$, infer the rule:
\begin{equation}
\label{DmBL:VCU:1}
\mbox{Being proved }\phi\rightarrow\psi,\mbox{ it is proved }\phi\counterfact\psi\;.
\end{equation}
It appears that Ax.2, Ax.4, Ax.5, Ax.6 and CR are recovered in DmBL, Ax.1 is weakened in DmBL and Ax.3 is not implemented in DmBL.
\\[5pt]
Conversely, b3 is implemented by VCU.
b2 is not implemented by VCU, but it could be shown that VCU completed by b4 implies b2.
b4 is not implemented by VCU.
b1 is obtained from (\ref{DmBL:VCU:1}), while weakened by $\Box\neg\phi$.
b5 is related to the notion of logical independence, which is not considered within VCU.
Then we have to point out three fundamental distinctions of DmBL compared to VCU:
\begin{enumerate}
\item\label{DmBL:Point:1} In DmBL, the negation commutes with the conditional (b4).
More generally, sub-universes are classical in DmBL,
\item\label{DmBL:Point:2} In DmBL, the deductions on the conditionals are often weakened by the hypothesis that \emph{the condition is not empty}; for example, $\Box\neg\phi$ in rule b1, or theorem $\vdash(\phi|\phi)\vee\Box\neg\phi$\,,
\item DmBL manipulates a notion of logical independence of the propositions.
\end{enumerate}
\emph{In fact, point~\ref{DmBL:Point:1} (commutation of the negation) makes point~\ref{DmBL:Point:2} (deduction weakened by the non-empty condition hypothesis) necessary.}
For example, $\bot\counterfact\top$ is derived from (\ref{DmBL:VCU:1});
by using both Ax.1 and the negation commutation, it is then deduced:
\begin{equation} \label{DmBL:VCU:eq1}
\top\equiv\bot\counterfact\bot\equiv\bot\counterfact\neg\top\equiv\neg(\bot\counterfact\top)\equiv\neg\top\equiv\bot\;,
\end{equation}
which is impossible.
Notice that this deduction is also done in DmBL, if we replace the ``weakened'' theorem $\vdash\Box\neg\phi\vee\Box(\phi|\phi)$ by the ``strong'' theorem $\vdash(\phi|\phi)$\,.
\\[5pt]
This example, based on VCU and DmBL, illustrates a fundamental difference between DmBL and other conditional logics.
DmBL considers $\bot$ as a singularity, and will be cautious with this case when inferring conditionals.
This principle is not just a logical artifact.
In fact, it is also deeply related to the notion of logical independence, as it appears in the proof of theorem \emph{Independence and proof}.
\section{Models}
\label{Section:Model:DmBL}
\subsection{Toward a Model}
In this paragraph, it is discussed about the link between Kripke models~\cite{blackburn,giordano1,giordano2} for DmBL/DmBL$_\ast$ and a more basic structure called \emph{conditional models}.\\[5pt]
\emph{From now on, $\mathcal{P}(W)$ denotes the set of all subsets of set $W$.}
\paragraph{Definition.}
A Kripke model for DmBL (respectively DmBL$_\ast$) is a quadruplet $(W,R,H,f)$, where $W$ is a set of worlds, $R\subset W\times W$ is an accessibility relation, $H:\mathcal{L}\longrightarrow\mathcal{P}(W)$ is an assignment function, $f:H(\mathcal{L})\times H(\mathcal{L})\longrightarrow\mathcal{P}(W)$ is a conditioning function, and verifying:
\begin{itemize}
\item $H(\neg\phi)=W\setminus H(\phi)$
and $H(\phi\rightarrow\psi)=\bigl(W\setminus H(\phi)\bigr)\cup H(\psi)$\,,
\item $H(\Box\phi)=\bigl\{t\in W\big/ \forall u\in W,\,(t,u)\in R\Rightarrow u\in H(\phi)\bigr\}$\,,
\item $H\bigl((\psi|\phi)\bigr)=f\bigl(H(\psi),H(\phi)\bigr)$\,,
\item $H(\phi)=W$ for any $\phi$ such that $\vdash\phi$ is an axiom of the form m3, b1, b2, b3, b4,\\
or b5 (respectively b5.weak.A and b5.weak.B).
\end{itemize}
It is noticed that the rules and axioms c$\ast$, modus ponens, m1 and m2 and are compliant with the model by construction.
\paragraph{Definition 2.} A conditional model for DmBL (respectively DmBL$_\ast$) is a quadruplet $(W,M,h,f)$ such that $W$ is a set of worlds, $M\subset\mathcal{P}(W)$ is a set of admissible propositions, $h:\Theta\longrightarrow M$ is an assignment function, $f:M\times M\longrightarrow M$ is a conditioning function, and verifying:
\begin{description}
\item[$\rien\quad\bullet$] $M$ is a Boolean sub-algebra of $\mathcal{P}(W)$, \emph{i.e.} $A\cap B\in M$ and $W\setminus A\in M$ for any $A,B\in M$\,,
\item[$\rien\quad\beta1$.] $A\subset B$ and $A\ne\emptyset$ imply $f(B,A)=W$\,, for any $A,B\in M$\,,
\item[$\rien\quad\beta2$.] $f(B\cup C,A)\subset f(B,A)\cup f(C,A)$\,, for any $A,B,C\in M$\,,
\item[$\rien\quad\beta3$.] $A\cap f(B,A) \subset B$\,, for any $A,B\in M$\,,
\item[$\rien\quad\beta4$.] $f(W\setminus B,A)=W\setminus f(B,A)$\,,
for any $A,B\in M$\,,
\item[$\rien\quad\beta5$.] $f(B,A)=B$ implies $f(A,B)=A$\\
(respectively $\beta5w.$ $f(B,A)=B$ implies  $f(B,W\setminus A)=B$)\,,
for any $A,B\in M$)\,.
\vspace{5pt}\end{description}
\emph{Remarks.} A conditional model does not implement the modalities.
In both models, the function $f$ is the representation of the conditional $(|)$\,.
\paragraph{Model transfer.}
Let $(W,M,h,f)$ be a conditional model for DmBL (respectively DmBL$_\ast$).
Let $R=W\times W$ and define $H$ by:
\begin{itemize}
\item $H(\theta)=h(\theta)$ for any $\theta\in\Theta$\,,
\item $H(\neg\phi)=W\setminus H(\phi)$ and $H(\phi\rightarrow\psi)=\bigl(W\setminus H(\phi)\bigr)\cup H(\psi)$\,,
\item $H(\phi)=W\Rightarrow H(\Box\phi)=W$ and $H(\phi)\ne W\Rightarrow H(\Box\phi)=\emptyset$\,,
\item $H\bigl((\psi|\phi)\bigr)=f\bigl(H(\psi),H(\phi)\bigr)$\,.
\end{itemize}
Then $(W,R,H,f)$ is a Kripke model for DmBL (respectively DmBL$_\ast$).
\\[5pt]
\emph{The proof is easy, but tedious.
It is detailed in appendix~\ref{proof:modtrans}.}
\\[5pt]
Conditional models are defined from the classical and conditional operators only.
In fact, such models have been set first for a non-modal construction of the Bayesian logic~\cite{Dambreville:DBL}.
In this paper a free conditional model is constructed for DmBL$_\ast$\,, with completeness results.
The conditional model is translated into a DmBL$_\ast$ Kripke model.
The derived model is of course not complete for DmBL$_\ast$ in regards to the modalities, but the completeness still holds in regards to the conditional operator.
The \emph{model transfer property} also suggests that the conditional operator is \emph{not} constructed from the modal operator: it is even possible to construct $(|)$ when $\Box$ is trivial in the model ($R=W\times W$ implies $H(\phi)=W$ or $H(\Box\phi)=\emptyset$)\,.
\\[5pt]
In section~\ref{Section:Proba:DmBL}, a model for DmBL is also derived but not constructed.
This model of DmBL is non-trivial (it is possible to extend any unconditioned probability over this model), but no completeness result is provided.
\subsection{Construction of a free conditional model for DmBL$_\ast$}
In the sequel, $\Theta$ is assumed to be finite.
A free conditional model for DmBL$_\ast$ will be constructed as a limit of partial models.
These models are constructed recursively, based on the iteration of $(|)$ on any propositions.
\subsubsection{Definition of partial models}
In this section are constructed a sequence $(\Omega_n, M_n, h_n, f_n,\Lambda_n)_{n\in\Nset}$ and a sequence of one-to-one morphisms $(\mu_n)_{n\in\Nset}$ such that:
\begin{itemize}
\item $M_n$ is a Boolean sub-algebra of $\mathcal{P}(\Omega_n)$, $h_n:\Theta\rightarrow M_n$ and $f_n:M_n\times M_n\rightarrow M_n$\,,
($f_n$ will be partially defined)
\item $\mu_n:M_n\rightarrow M_{n+1}$ is such that $\mu_n(A\cap B)=\mu_n(A)\cap\mu_n(B)$, $\mu_n(\Omega_n\setminus A)=\Omega_{n+1}\setminus\mu_n(A)$
and $\forall\theta\in\Theta,\,h_{n+1}(\theta)=\mu_n\bigl(h_n(\theta)\bigr)$\,,
(as a consequence, $\mu_n$ is a Boolean morphism)
\item For any $A,B\in M_n$ such that $f_n(B,A)$ is defined, then $f_{n+1}\bigl(\mu_n(B),\mu_n(A)\bigr)$ is defined and $f_{n+1}\bigl(\mu_n(B),\mu_n(A)\bigr)=\mu_n\bigl(f_n(B,A)\bigr)$\,,
\item $\Lambda_n$ is a list of elements of $M_n$, which is used as a task list of the construction (refer to the subsequent paragraphs).
\end{itemize}
\emph{Remark.}
The functions $f_n$ represent the partial construction of $(|)$\,.
The morphisms $\mu_n$ characterize the ``inclusion'' of the partial models.
\\[5pt]
In a subsequent section, a conditional model $\bigl(\Omega_\infty,M_\infty,h_\infty,f_\infty\bigr)$ will be defined as the limit of $(\Omega_n, M_n, h_n, f_n)_{n\in\Nset}$ associated to $(\mu_n)_{n\in\Nset}$\,.
\paragraph{Notations and definitions.}
For any $A\in M_n$, it is defined $\sim A=\Omega_n\setminus A$\,.\\[5pt]
Any singleton $\{\omega\}$ may be denoted $\omega$ if the context is not ambiguous.\\[5pt]
For any $m>n$ and $A\in M_n$\,, it is defined $A_{[m]}=\mu_{m-1}\circ\dots\circ\mu_{n}(A)$\,.\\[5pt]
The Cartesian product of sets $A$ and $B$ is denoted $A\times B$\,;
the functions $\mathrm{id}$ and $T$ are defined over pairs by $\mathrm{id}(x,y)=(x,y)$ and $T(x,y)=(y,x)$\,;
for a set of pairs $C$, the abbreviation $(\mathrm{id}\cup T)(C)=\mathrm{id}(C)\cup T(C)$ is also used.
\paragraph{Initialization.}
Define $(\Omega_0, M_0, h_0, f_0,\Lambda_0)$ by:
\begin{itemize}
\item $\Omega_0=\{0,1\}^{\Theta}$,
\item $M_0=\mathcal{P}(\Omega_0)$,
\item $\forall\theta\in\Theta,\,h_0(\theta)=\bigl\{(\delta_\tau)_{\tau\in\Theta}\in\Omega_0\,\big/\,\delta_\theta=1\bigr\}$\,,
\item $f_0(A,\emptyset)=f_0(A,\Omega_0)=A$ for any $A\in M_0$\,,
\item $\Lambda_0=(s_0,f_0,\lambda_0)$ is a list defined by $s_0=0$, $f_0=\mathrm{card}(M_0)-2$ and $\lambda_0$ is a one-to-one mapping from $[\![s_0,f_0-1]\!]$ to $M_0\setminus\{\emptyset,\Omega_0\}$\,,
such that $\lambda_0(2t)=\sim\lambda_0(2t+1)\,,\;\forall t$\,.
\end{itemize}
\paragraph{Step $n$ to step $n+1$.}\rien\\
Let $(\Omega_k, M_k, h_k, f_k,\Lambda_k)_{0\le k\le n}$ and the one-to-one morphisms $(\mu_k)_{0\le k\le n-1}$ be constructed.\\[5pt]
Notice that $\lambda_n(s_n)=\sim\lambda_n(s_n+1)$ by construction of $\Lambda_n$\,.\\
\underline{Define $b_{n}=\lambda_n(s_n)$}\,.\\
Then, construct the set $I_n$ and the sequences $\Gamma_n(i),\Pi_n(i)|_{i\in I_n}$ according to the cases:
\subparagraph{Case 0.} There is $m<n$ such that $\{b_{m[n]},\sim b_{m[n]}\}=\{b_n,\sim b_n\}$\,.\\
Then, notice that $b_n=b_{m[n]}$ by the subsequent construction of $\Lambda$\,.\\
Let $\nu$ be the greatest of such $m$;
then define $I_n=\mu_\nu(b_{\nu})\times\sim \mu_\nu(b_{\nu})$\,,\\
$\Pi_n(\omega,\omega')=f_n(\omega'_{[n]},\sim b_n)\cap\omega_{[n]}$ and $\Gamma_n(\omega,\omega')=f_n(\omega_{[n]},b_n)\cap\omega'_{[n]}$ for any $(\omega,\omega')\in I_n$\,.$^\dagger$
\\[5pt]
$^\dagger$Remark: case 0 means that the construction of $f(\cdot,b_n)$ and of $f(\cdot,\sim b_n)$ has already begun over the propositions of $M_{\nu+1}$.
\subparagraph{Case 1.} Case 0 does not hold;\\
Define $I_n=\{b_n\}$\,, $\Pi_n(i)=i$ 
and $\Gamma_n(i)=\sim i$ for any $i\in I_n$\,.
\\[5pt]
Remark: case 1 means that $f(\cdot,b_n)$ and $f(\cdot,\sim b_n)$ are constructed for the first time.
\subparagraph{Setting.}
$(\Omega_{n+1}, M_{n+1}, h_{n+1}, f_{n+1},\Lambda_{n+1})$ and $\mu_n$ are defined by:
\begin{itemize}
\item $\mu_n(A)=\bigcup_{i\in I_n} \biggl(\Bigl(\bigl(A\cap\Pi_n(i)\bigr)\times\Gamma_n(i)\Bigr)
\cup \Bigl(\bigl(A\cap\Gamma_n(i)\bigr)\times\Pi_n(i)\Bigr)\biggr)$ for any $A\in M_n$\,,
\item $\Omega_{n+1}=\mu_n(\Omega_n)$\,,
\item $\forall\theta\in\Theta,\,h_{n+1}(\theta)=\mu_n\bigl(h_n(\theta)\bigr)$\,,
\item $M_{n+1}=\mathcal{P}(\Omega_{n+1})$\,,
\item $f_{n+1}(A,\emptyset)=f_{n+1}(A,\Omega_{n+1})=A$ for any $A\in M_{n+1}$\,,
\item For any $A\in M_n\setminus\{b_n,\sim b_n,\emptyset,\Omega_n\}$ and any $B\in M_n$ such that $f_n(B,A)$ is defined, then $f_{n+1}\bigl(\mu_n(B),\mu_n(A)\bigr)$ is defined and $f_{n+1}\bigl(\mu_n(B),\mu_n(A)\bigr)=\mu_n\bigl(f_n(B,A)\bigr)$\,,
\item For any $A\in M_{n+1}$\,, set
$
f_{n+1}\bigl(A,\mu_n(b_n)\bigr)=(\mathrm{id}\cup T)
\biggl(A\cap\Bigl(\bigcup_{i\in I_n}\bigl(\Pi_n(i)\times\Gamma_n(i)\bigr)\Bigr)\biggr)
$\\[5pt]
and
$
f_{n+1}\bigl(A,\sim\mu_n(b_n)\bigr)=(\mathrm{id}\cup T)
\biggl(A\cap\Bigl(\bigcup_{i\in I_n}\bigl(\Gamma_n(i)\times\Pi_n(i)\bigr)\Bigr)\biggr)
\;,
$
\item $\Lambda_{n+1}=(s_{n+1},f_{n+1},\lambda_{n+1})$ is such that:
\begin{itemize}
\item $s_{n+1}=s_n+2=2n+2$, $f_{n+1}=s_{n+1}+\mathrm{card}(M_{n+1})-2$\,,
\item $\lambda_{n+1}$ is a one-to-one mapping from $[\![s_{n+1},f_{n+1}-1]\!]$ to $M_{n+1}\setminus\{\emptyset,\Omega_{n+1}\}$\,,
\item $\lambda_{n+1}(t)=\mu_n\bigl(\lambda_n(t)\bigr)$ for any $t\in[\![s_{n+1},f_{n}-1]\!]$\,,
\item $\lambda_{n+1}(2t)=\sim\lambda_{n+1}(2t+1)\,,\;\forall t$
and $\lambda_{n+1}(f_{n+1}-2)=\mu_n(b_n)$\,.
\end{itemize}
\emph{This definition ensures a cyclic and full construction of $f(\cdot,b_n)$ and $f(\cdot,\sim b_n)$}.\footnote{In regards to the mapping $\mu$, the list $\Lambda_{n+1}$ is the list $\Lambda_n$ plus any propositions of $M_{n+1}\setminus\{\emptyset,\Omega_{n+1}\}$ which are not/\emph{no more} listed in $\Lambda_n$\,.}
\end{itemize}
The first steps of the model construction are illustrated by a simple example in appendix~\ref{BayesMod:Comp&Ex}.
%
\paragraph{Short explanation of the model.}
In fact, $(\omega,\omega')\in\Pi_n(i)\times\Gamma_n(i)$ should be interpreted as $\omega\wedge(\omega'|\neg b_n)$, while $(\omega',\omega)\in\Gamma_n(i)\times\Pi_n(i)$ should be interpreted as $\omega'\wedge(\omega|b_n)$.
The reader should compare this construction to the proof of completeness in appendix~\ref{Appendix:ProofOfAlmostCompletude} for a better comprehension of the mechanisms of the model.
\subsubsection{Properties of $(\Omega_n, M_n, h_n, f_n,\Lambda_n,\mu_n)_{n\in\Nset}$}
\label{Omega:properties}
It is proved recursively:
\begin{description}
\item[$\rien\quad\bullet_\mu$] $\mu_n:M_n\rightarrow M_{n+1}$ is a one-to-one Boolean morphism,
\item[$\rien\quad\bullet_f$] If $A,B\in M_n$ and $f_n(B,A)$ is defined, then $f_{n+1}\bigl(\mu_n(B),\mu_n(A)\bigr)=\mu_n\bigl(f_n(B,A)\bigr)$\,,
\item[$\rien\quad\tilde\beta1$.] Let $A,B\in M_n$ such that $f_n(B,A)$ is defined.
\\Then $A\subset B$ and $A\ne\emptyset$ imply $f_n(B,A)=\Omega_n$\,,
\item[$\rien\quad\tilde\beta2$.] Let $A,B,C\in M_n$ such that $f_n(B,A)$, $f_n(C,A)$ and $f_n(B\cup C,A)$ are defined.
\\Then $f_n(B\cup C,A)= f_n(B,A)\cup f_n(C,A)$\,,
\item[$\rien\quad\tilde\beta3$.] Let $A,B\in M_n$ such that $f_n(B,A)$ is defined.
\\Then $A\cap f_n(B,A) = A\cap B$\,,
\item[$\rien\quad\tilde\beta4$.] Let $A,B\in M_n$ such that $f_n(B,A)$ and $f_n(\sim B,A)$ are defined.
\\Then $f_n(\sim B,A)=\sim f_n(B,A)$\,,
\item[$\rien\quad\tilde\beta5w$.] Let $A,B\in M_n$ such that $f_n(B,A)$ and $f_n(B,\sim A)$ are defined.
\\Then $f_n(B,A)=B$ implies $f_n(B,\sim A)=B$\,.
\end{description}
\emph{Proofs are given in appendix~\ref{Appendix:MainProof}.}
\subsubsection{Limit}
The limit $\bigl(\Omega_\infty,M_\infty,h_\infty,f_\infty\bigr)$ is defined as follows:
\begin{itemize}
\item Set $\displaystyle\Omega_\infty=\left\{\left.(\omega_n|_{n\in\Nset})\in\prod_{n\in\Nset}\Omega_n\;\right/\;\forall n\in\Nset,\,\omega_{n+1}\in\mu_n(\omega_n)\right\}$\,;\\[10pt]
\emph{Useful definitions:}
\begin{itemize}
\item For any $n\in\Nset$ and any $A\in M_n$\,, $\displaystyle A_\infty=\bigl\{(\omega_k|_{k\in\Nset})\in\Omega_\infty
\;\big/\;\omega_n\in A\bigr\}$\,.
\emph{
The subset $A_\infty$ is a mapping of $A$ within $\Omega_\infty$\,.
It is noticed that this mapping is invariant with $\mu$, \emph{i.e.} $(A_{[m]})_\infty=A_\infty$ for $m>n$\,,
}
\item For any $n\in\Nset$\,, $\displaystyle M_{n:\infty}=\bigl\{A_\infty\;\big/\;A\in M_n\bigr\}$\,.
\emph{
The structure $M_{n:\infty}$ is an isomorphic mapping of the structure $M_n$ within $\Omega_\infty$\,.
It is noticed that $(M_{n:\infty}|_{n\in\Nset})$ is a monotonic sequence, \emph{i.e.} $M_{n:\infty}\subset M_{n+1:\infty}$\,,
}
\end{itemize}
\item Set $\displaystyle M_\infty=\bigcup_{n\in\Nset}M_{n:\infty}$\,,
\item Set $h_\infty(\theta)=\bigl(h_0(\theta)\bigr)_\infty$ for any $\theta\in\Theta$\,,
\item Let $A,B\in M_\infty$\,.
Then there is $n$ and $a,b\in M_n$ such that $A=a_\infty$\,, $B=b_\infty$ and $f_n(b,a)$ is defined (subsequent proposition).
Set $f_\infty(B,A)=\bigl(f_{n}(b,a)\bigr)_\infty$\,.
\end{itemize}
This definition is justified by the following propositions:
\paragraph{Proposition 1.}
For any $n\in\Nset$, $m>n$ and $A\in M_n$, $(A_{[m]})_\infty=A_\infty$\,.
\begin{description}
\item[Proof.] By definition of $\Omega_\infty$,
$$
(A_{[m]})_\infty=\bigl\{(\omega_k|_{k\in\Nset})\in\Omega_\infty
\;\big/\;\omega_m\in A_{[m]}\bigr\}=\bigl\{(\omega_k|_{k\in\Nset})\in\Omega_\infty
\;\big/\;\omega_n\in A\bigr\}
=A_\infty\;.
$$
\item[$\Box\Box\Box$]\rien
\end{description}
\emph{Corollary} $M_{n:\infty}\subset M_{n+1:\infty}$\,.
\paragraph{Proposition 2.} $M_{n:\infty}$ is a Boolean subalgebra of $\mathcal{P}(\Omega_\infty)$ and is isomorph to $M_n$ by the morphism $A\mapsto A_\infty$. As a consequence, $M_{\infty}$ is a Boolean subalgebra of $\mathcal{P}(\Omega_\infty)$\,.
\\[5pt]
\emph{Proof is obvious from the definition of $M_{n:\infty}$}\,.\\[5pt]
From now on, $M_n$ will be considered as a subalgebra of $\mathcal{P}(\Omega_\infty)$.
\paragraph{Proposition 3.} Let $A,B\in M_\infty$\,.
Then there is $n$ and $a,b\in M_n$ such that $A=a_\infty$\,, $B=b_\infty$ and $f_n(b,a)$ is defined.
\begin{description}
\item[Proof.]
Since $(M_{k:\infty}|_{k\in\Nset})$ is a monotonic sequence\,, there is $m\in\Nset$ such that $A,B\in M_{m:\infty}$\,.\\
Let $a,b\in M_m$ be such that $A=a_\infty$ and $B=b_\infty$\,.\\
By definition of the list $\Lambda_m$, there is $p\in[\![s_{m},f_{m}-1]\!]$ such that $\lambda_p=a$\,.\\
As a consequence, $f_{p+1}(b_{[p+1]},a_{[p+1]})$ exists.\\
But then hold $A=(a_{[p+1]})_\infty$ and $B=(b_{[p+1]})_\infty$\,.\\
Finally $n=p+1$ answers to the proposition.
\item[$\Box\Box\Box$]\rien
\end{description}
\paragraph{Proposition 4.} The definition of $f_\infty$ does not depend on the choice of $n$.
\begin{description}
\item[Proof.]
Let $n$, $m>n$, $a,b\in M_n$ and $c,d\in M_m$ such that $A=a_\infty=c_\infty$ and $B=b_\infty=d_\infty$\,.\\
Assume also that $f_n(b,a)$ and $f_m(d,c)$ exist.\\[3pt]
Then $(a_{[m]})_{\infty}=c_\infty$ and $(b_{[m]})_{\infty}=d_\infty$.\\
Since $M_m$ and $M_{m:\infty}$ are isomorph, it follows $a_{[m]}=c$ and $b_{[m]}=d$\,.\\
But it is derived from $\bullet_f$ that $f_m(b_{[m]},a_{[m]})=f_n(b,a)_{[m]}$\,.\\
Finally $\bigl(f_n(b,a)\bigr)_\infty=\bigl(f_m(d,c)\bigr)_\infty$\,.
\item[$\Box\Box\Box$]\rien
\end{description}
\paragraph{Proposition 5.} $\bigl(\Omega_\infty,M_\infty,h_\infty,f_\infty\bigr)$  verifies $\beta1$, $\beta2$, $\beta3$, $\beta4$ and $\beta5w$.\\[5pt]
The properties $\beta\ast$ are inherited from $M_n,f_n|_{n\in\Nset}$\,, by means of the properties $\tilde\beta\ast$.
\paragraph{Conclusion.} \underline{$\bigl(\Omega_\infty,M_\infty,h_\infty,f_\infty\bigr)$ is a conditional model for DmBL$_\ast$.}
\subsubsection{Implied Kripke model for DmBL$_\ast$}
By means of the \emph{Model transfer} property, a Kripke model for DmBL$_\ast$ is derived from $\bigl(\Omega_\infty,M_\infty,h_\infty,f_\infty\bigr)$.
\underline{This Kripke model is denoted $\mathcal{B}=\bigl(\Omega_\infty,R_{\mathcal{B}},H_{\mathcal{B}},f_\infty\bigr)$\,.}
\subsubsection{Completeness for the conditional operator}
It is above the scope of this work to construct a model of DmBL$_\ast$, which is complete for both the modal and the conditional operators.
However, it is shown here that $\mathcal{B}$ is a model of DmBL$_\ast$, which is complete for the conditional operator.
\paragraph{Proposition 1.}
By construction, $(\Omega_\infty,H_{\mathcal{B}})$ is a complete model for the classical logic $C$, when $H_{\mathcal{B}}$ is restricted to the propositions of $\mathcal{L}_C$.
\paragraph{Proposition 2.}
Let $\phi\in\mathcal{L}$ be a proposition constructed \underline{without $\Box$ or $\Diamond$}\,.
Then $\vdash\phi$ in DmBL$_\ast$ if and only if $H_{\mathcal{B}}(\phi)=\Omega_\infty$\,.
\\[5pt]
\emph{Proof is done in appendix~\ref{Appendix:ProofOfAlmostCompletude}\,.}
\\[5pt]
Proposition 2 expresses that $\mathcal{B}$ is complete for the conditional operator.
\subsection{Coherence properties}
The model $\mathcal{B}$ clearly shows that DmBL$_\ast$ is coherent.
It also demonstrates that the conditional operator $(|)$ is not trivial.
Since $(\Omega_\infty,H_{\mathcal{B}})$ is a complete model for $C$, DmBL$_\ast$ is an extension of the classical logic:
$\vdash\phi$ implies $\vdash_C\phi$, for any $\phi\in\mathcal{L}_C$.
But a stronger property holds:
\paragraph{Non-distortion.}
Let $\phi\in\mathcal{L}_C$.
Assume that $\vdash\Box\phi\vee\Box\neg\phi$ in DmBL$_\ast$.
Then $\vdash_C\phi$ or $\vdash_C\neg\phi$.
\\[5pt]
\emph{Interpretation: DmBL$_\ast$ does not ``distort'' the \emph{classical} propositions.
More precisely, a property like $\vdash\Box\phi\vee\Box\neg\phi$ would add some knowledge about $\phi$, since it says that either $\phi$ or $\neg\phi$ is ``sure''.
But the \emph{non-distortion} just tells that such property is impossible unless there is a trivial knowledge about $\phi$ within the classical logic.}
\begin{description}
\item[Proof.]
Assume $\vdash\Box\phi\vee\Box\neg\phi$\,.\\
Since $\mathcal{B}$ is a model for DmBL$_\ast$, it comes $H_{\mathcal{B}}(\Box\phi\vee\Box\neg\phi)=\Omega_\infty$\,.\\
Then $H_{\mathcal{B}}(\Box\phi)\cup H_{\mathcal{B}}(\Box\neg\phi)=\Omega_\infty$\,.\\
Since $H_{\mathcal{B}}(\Box \psi)=\emptyset$ or $\Omega_\infty$ for any $\psi\in\mathcal{L}$\,, it comes $H_{\mathcal{B}}(\Box\phi)=\Omega_\infty$ or $H_{\mathcal{B}}(\Box\neg\phi)=\Omega_\infty$\,.\\
At last, $H_{\mathcal{B}}(\phi)=\Omega_\infty$ or $H_{\mathcal{B}}(\neg\phi)=\Omega_\infty$\,.\\
But $(\Omega_\infty,H_{\mathcal{B}})$ is a complete Boolean model for $C$\,, which implies $\vdash_C\phi$ or $\vdash_C\neg\phi$\,.
\item[$\Box\Box\Box$]\rien
\end{description}
Another \emph{non-distortion} property is derived subsequently in the context of probabilistic DmBL$_\ast$.
\section{Extension of probability}
\label{Section:Proba:DmBL}
\subsection{Probability over propositions,} a minimal$^\dagger$ definition.\\
{\small $\dagger$ This definition is related to finite probabilities and excludes any Bayesian consideration.}
\\[7pt]
Probabilities are classically defined over measurable sets.
However, this is only a manner to model the notion of probability, which is essentially an additive measure of the belief of logical propositions \cite{paass}.
Probability could be defined without reference to the measure theory, at least when the propositions are countable.
The notion of probability is explained now within a strict propositional formalism.
Conditional probabilities are excluded from this definition, but the notion of independence is considered.
\vspace{5pt}\\
Intuitively, a probability over a space of logical propositions is a measure of belief which is additive (disjoint propositions are adding their chances) and increasing with the propositions.
This measure should be zeroed for the contradiction and set to $1$ for the tautology.
Moreover, \emph{a probability is a multiplicative measure for independent propositions}.
%
\paragraph{Definition for classical propositions.}
A probability $\pi$ over $C$ is a $\Rset^+$ valued function such that for any propositions $\phi$ and $\psi$ of $\mathcal{L}_C$\,:
\begin{description}
\item[\rien$\quad$\emph{Equivalence.}]$\phi\equiv_C\psi$ implies $\pi(\phi)=\pi(\psi)$\,,
\item[\rien$\quad$\emph{Additivity.}]$\pi(\phi\wedge\psi)+\pi(\phi\vee\psi)=\pi(\phi)+\pi(\psi)$\,,
\item[\rien$\quad$\emph{Coherence.}]$\pi(\bot)=0$\,,
\item[\rien$\quad$\emph{Finiteness.}]$\pi(\top)=1$\,.
\end{description}
\subparagraph{Property.}
The coherence and additivity imply the increase of $\pi$:
\begin{description}
\item[\rien$\quad$\emph{Increase.}]$\pi(\phi\wedge\psi)\le \pi(\phi)$\,.
\end{description}
\begin{description}
\item[Proof.] 
Since $\phi\equiv_C(\phi\wedge\psi)\vee(\phi\wedge\neg\psi)$ and $(\phi\wedge\psi)\wedge(\phi\wedge\neg\psi)\equiv_C\bot$, the additivity implies:
$$
\pi(\phi)+\pi(\bot)=\pi(\phi\wedge\psi)+\pi(\phi\wedge\neg\psi)\;.
$$
From the coherence $\pi(\bot)=0$\,,
it is deduced $\pi(\phi)=\pi(\phi\wedge\psi)+\pi(\phi\wedge\neg\psi)$\,.\\
Since $\pi$ is non-negatively valued, $\pi(\phi)\ge \pi(\phi\wedge\psi)$\,.
\item[$\Box\Box\Box$]\rien
\end{description}
\paragraph{Definition for DmBL/DmBL$_\ast$.}
In this case, we have to deal with independence notions.\\[5pt]
A probability $P$ over DmBL/DmBL$_\ast$ is a $\Rset^+$ valued function, which verifies (replace $\equiv_C$ by $\equiv$ and $\pi$ by $P$) \emph{equivalence}, \emph{additivity}, \emph{coherence}, \emph{finiteness} and:
\begin{description}
\item[\rien$\quad$\emph{Multiplicativity.}]$\vdash\phi\times\psi$ implies $P(\phi\wedge\psi)=P(\phi)P(\psi)$\,.
\end{description}
for any propositions $\phi$ and $\psi$ of $\mathcal{L}$\,. 
\subsection{Probability extension over DmBL$_\ast$}
\label{ProbExt:wDmBL}
\paragraph{Property.}
Let $\pi$ be a probability defined over $C$\,, the classical logic, such that $\pi(\phi)>0$ for any $\phi\not\equiv_C\bot$.
Then, there is a (multiplicative) probability $\overline{\pi}$ defined over DmBL$_\ast$ such that $\overline{\pi}(\phi)=\pi(\phi)$ for any classical proposition $\phi\in\mathcal{L}_C$\,.
\\\\
\emph{Remark: this is another non-distortion property, since the construction of DmBL$_\ast$ puts no constraint over probabilistic classical propositions.}
\\[5pt]
Proof is done in appendix~\ref{Appendix:Probabilition}.
\paragraph{Corollary.}
Let $\pi$ be a probability defined over $C$\,.
Then, there is a (multiplicative) probability $\overline{\pi}$ defined over DmBL$_\ast$ such that $\overline{\pi}(\phi)=\pi(\phi)$ for any $\phi\in\mathcal{L}_C$\,.
\begin{description}
\item[Proof.]
Let $\Sigma=\left\{\left.\bigwedge_{\theta\in\Theta}\epsilon_\theta\;\right/\;\epsilon\in\prod_{\theta\in\Theta}\{\theta,\neg\theta\}\right\}$\,.\\
For any real number $e >0$\,, define the probability $\pi_e $ over $\mathcal{L}_C$ by:
$$
\forall\sigma\in\Sigma\,,\; \pi_e (\sigma)=\frac{e }{\mathrm{card}(\Sigma)}+(1-e )\pi(\sigma)
\;. 
$$
Let $\overline{\pi_e}$ be the extension of  $\pi_e$ over DmBL$_\ast$ as constructed in appendix~\ref{Appendix:Probabilition}.\\
By~\ref{AppC:Conclude}\,, there is a rational function $R_\phi$ such that $\overline{\pi_e}(\phi)=R_\phi(e )$ for any $\phi\in\mathcal{L}$\,.\\
Now $0\le R_\phi(e )\le 1$\,;
since $R_\phi(e )$ is rational and bounded, $\lim_{e \rightarrow 0+}R_\phi(e )$ exists.\\
Define $\overline{\pi}(\phi)=\lim_{e \rightarrow 0+}R_\phi(e )$\,, for any $\phi\in\mathcal{L}$.\\
The additivity, coherence, finiteness and multiplicativity are obviously inherited by $\overline{\pi}$.\\
At last, it is clear that $\overline{\pi}(\sigma)=\pi(\sigma)$ for any $\sigma\in\Sigma$\,.
\item[$\Box\Box\Box$]\rien
\end{description}
\subsection{Model and probability extension for DmBL}
Let $\mathcal{K}$ be the set of all (multiplicative) probabilities $P$ over DmBL$_\ast$ such that $P(\phi)>0$ for any $\phi\not\equiv\bot$\,,
and define the sequences $\mathcal{K}(\phi)=(P(\phi))_{P\in\mathcal{K}}$ for any $\phi\in\mathcal{L}$\,.\\
Then define
$\mathcal{L}_{\mathcal{K}}=\mathcal{K}(\mathcal{L})=\bigl\{\mathcal{K}(\phi)\;\big/\;\phi\in\mathcal{L}\bigr\}\;;$
The space $\mathcal{L}_{\mathcal{K}}$ is thus a subset of ${\Rset^+}^{\mathcal{K}}$\,.\\
The operators $\neg$, $\wedge$ and $(|)$ are canonically implied over $\mathcal{L}_{\mathcal{K}}$\,:\footnote{The operators $\vee$ and $\rightarrow$ are derived from $\wedge$ and $\neg$ as usually; modalities are not considered.}
$$
\neg\mathcal{K}(\phi)=\mathcal{K}(\neg\phi)\,,\ 
\mathcal{K}(\phi)\wedge\mathcal{K}(\psi)=\mathcal{K}(\phi\wedge\psi)
\quad\mbox{and}\quad
\bigl(\mathcal{K}(\psi)\big|\mathcal{K}(\phi)\bigr)=
\mathcal{K}\bigl((\psi|\phi)\bigr)\;.
$$
Since any $P\in\mathcal{K}$ verifies the equivalence property, it comes $\mathcal{K}(\phi)=\mathcal{K}(\psi)$ when $\phi\equiv\psi$ in DmBL$_\ast$.
As a direct consequence, $\bigl(\mathcal{L}_{\mathcal{K}},\neg,\wedge,(|)\bigr)$ is a conditional-like model of DmBL$_\ast$ (the structure is a Boolean algebra but not derived from set operators. 
This is the only difference with conditional models).
\paragraph{Property.}
$\bigl(\mathcal{L}_{\mathcal{K}},\neg,\wedge,(|)\bigr)$ is a conditional-like model of DmBL.
\begin{description}
\item[Proof.]Let $P\in\mathcal{K}$\,; $P$ is multiplicative.\\
Since $\vdash(\psi|\phi)\times\phi$ and $(\psi|\phi)\wedge\phi\equiv\psi\wedge\phi$ in DmBL$_\ast$, it comes $P\bigl((\psi|\phi)\bigr)P(\phi)=P(\psi\wedge\phi)$\,.
\\[3pt]
Now assume $\bigl(\mathcal{K}(\psi)\big|\mathcal{K}(\phi)\bigr)=\mathcal{K}(\psi)$\,, with $\psi\not\equiv\bot$\,.\\
Then $\mathcal{K}\bigl((\psi|\phi)\bigr)=\mathcal{K}(\psi)$, and $P\bigl((\psi|\phi)\bigr)=P(\psi)$ for any $P\in\mathcal{K}$\,.\\
Then $P(\phi)=\frac{P(\psi\wedge\phi)}{P\bigl((\psi|\phi)\bigr)}=\frac{P(\psi\wedge\phi)}{P(\psi)}=P\bigl((\phi|\psi)\bigr)$ for any $P\in\mathcal{K}$\,,
\\
and $\bigl(\mathcal{K}(\phi)\big|\mathcal{K}(\psi)\bigr)=\mathcal{K}\bigl((\phi|\psi)\bigr)=\mathcal{K}(\phi)$\,.\\
Since moreover $(\phi|\bot)\equiv\phi$ and $(\bot|\phi)\equiv\bot$ in DmBL$_\ast$, the model verifies $\beta5$\,.
\item[$\Box\Box\Box$]\rien
\end{description}
Notice that it was only needed the equivalence and multiplicative properties for the elements of $\mathcal{K}$\,.
It is thus possible to construct a more general model by relaxing $\mathcal{K}$\,.
\paragraph{Probability extension.}
For any $\mathcal{K}(\phi)\in\mathcal{L}_{\mathcal{K}}$ and any $P\in\mathcal{K}$\,, define the $\Rset^+$-valued mapping ${\widehat P}\bigl(\mathcal{K}(\phi)\bigr)=P(\phi)$ (this mapping, a projection, is indeed well defined).\\
By construction, ${\widehat P}$ is naturally a multiplicative probability over $\mathcal{L}_{\mathcal{K}}$\,.
Moreover, the probability extensions defined in appendix~\ref{Appendix:Probabilition} are also elements of $\mathcal{K}$\,.
As a consequence, the deductions of section~\ref{ProbExt:wDmBL} are still working for $\mathcal{L}_{\mathcal{K}}$\,.
The extension property is thus derived:
\begin{quote}
Let $\pi$ be a probability defined over $C$\,.
Then, there is a (multiplicative) probability $\overline{\pi}$ defined over DmBL such that $\overline{\pi}(\phi)=\pi(\phi)$ for any $\phi\in\mathcal{L}_C$\,.
\end{quote}
\paragraph{Non-distortion.}
Let $\phi$ be a classical proposition.
Assume that $\vdash\Box\phi\vee\Box\neg\phi$ in DmBL.
Then $\vdash_C\phi$ or $\vdash_C\neg\phi$.
\begin{description}
\item[Proof.]
Consider the Kripke model for DmBL derived from the conditional model $\bigl(\mathcal{L}_{\mathcal{K}},\neg,\wedge,(|)\bigr)$\,.
\\
In this model, the value of $H(\Box\phi)$ is either $\mathcal{K}(\bot)$ or $\mathcal{K}(\top)$\,.\\
Then, $H(\Box\phi\vee\Box\neg\phi)=\mathcal{K}(\top)$ implies $H(\Box\phi)=\mathcal{K}(\top)$ or $H(\Box\neg\phi)=\mathcal{K}(\top)$\,.\\
Then $H(\phi)=\mathcal{K}(\top)$ or $H(\neg\phi)=\mathcal{K}(\top)$\,.\\
It follows $\forall P\in\mathcal{K}\,,\; P(\phi)=1$ or $\forall P\in\mathcal{K}\,,\; P(\neg\phi)=1$\,, and by the probability extension:
$\forall \pi\,,\; \pi(\phi)=1$ or $\forall \pi\,,\; \pi(\neg\phi)=1$\,, where $\pi$ denotes any probability over $C$\,.\\
At last, $\vdash_C\phi$ or $\vdash_C\neg\phi$\,.
\item[$\Box\Box\Box$]\rien
\end{description}
\subsection{Properties of the conditional}
\paragraph{Bayes inference.}
Assume a (multiplicative) probability $P$ defined over DmBL/DmBL$_\ast$.
Define $P(\psi|\phi)$ as an abbreviation for $P\bigl((\psi|\phi)\bigr)$\,.
Then:
$$
P(\psi|\phi)P(\phi)=P(\phi\wedge\psi)\;.
$$
\begin{description}
\item[Proof.]A consequence of \mbox{$(\psi|\phi)\wedge\phi\equiv\phi\wedge\psi$} and \mbox{$\vdash(\psi|\phi)\times\phi$}\,.
\item[$\Box\Box\Box$]\rien
\end{description}
As a corollary, it is also deduced
$P\bigl((\psi|\phi)\big|(\eta|\zeta)\bigr)P(\eta|\zeta)=P\bigl((\psi|\phi)\wedge(\eta|\zeta)\bigr)$\,.
It is recalled that the closure of CEA DGNW fails on this relation (refer to the introduction).
\paragraph{About Lewis' triviality.}
The previous extension theorems have shown that for any probability $\pi$ defined over $C$\,, it is possible to construct a (multiplicative) probability $\overline{\pi}$ over DmBL which extends $\pi$.
This result by itself shows that DmBL avoids Lewis' triviality.
But a deeper explanation seems necessary.
\\[5pt]
Assume $\phi\in\mathcal{L}_C$ and define the probability $\pi_\phi$ over $C$ by $\pi_\phi=\pi(\cdot|\phi)$\,.
Let $\overline{\pi_\phi}$ be the extension of $\pi_\phi$ over DmBL.
It happens that $\overline{\pi_\phi}\ne \overline{\pi}(\cdot|\phi)$\,, which implies that Lewis' triviality does not work anymore.
It is noticed that although $\overline{\pi}(\cdot|\phi)$ is a probability over DmBL in the classical meaning (it is \emph{additive}, \emph{coherent} and \emph{finite}), it is not necessarily \emph{multiplicative}.
\begin{quote}
Conditional probabilities do not maintain the logical independence and the conditioning.
\end{quote}
This limitation is unavoidable: otherwise the derivation~(\ref{Eq:DmBL:v2:Lewis:1}) of the triviality is possible, even if $(\psi|\phi)$ is not equivalent to a classical proposition.
\section{Conclusion}
\label{Fus2004::Sec:8}
In this contribution, the conditional logics DmBL and DmBL$_\ast$, a slight relaxation of DmBL, have been defined and studied.
These logics have been introduced as an abstraction and extrapolation of general probabilistic properties.
DmBL and DmBL$_\ast$ implement the essential ingredients of the Bayesian inference, including the classical nature to the sub-universe, the inference property and a related concept of logical independence.
For this reason, DmBL and DmBL$_\ast$ extend and refine the existing logical approaches of the Bayesian inference.
\\[5pt]
The logics are coherent and non-trivial.
A model has been constructed for the logic DmBL$_\ast$, which is complete in regards to the conditionals.
It has been shown that any probability over the classical propositions could be extended to DmBL/DmBL$_\ast$, in compliance with the independence relation.
Then, the probabilistic Bayesian rule has been recovered from DmBL/DmBL$_\ast$.
\\[5pt]
There are still many open questions.
For example, it is certainly possible to bring some enrichment to the conditional of DmBL, by means of additional axioms.
Is it possible to recover some specific equivalences of other existing systems?
From the strict logical viewpoint, the Deterministic modal Bayesian Logic offers also some interesting properties.
For example, the notion of independence in DmBL have nice logical consequences in the deductions (\emph{e.g.} regularity with an inference).
This property should be of interest in mathematical logic.
\small
%

%
\appendix
\section{Proof: the logical theorems}
\label{proof:logth}
%
%
\paragraph{Axioms order.}\rien\\
From b5, it is deduced $\psi\times\neg\phi\equiv\neg\phi\times\psi$\,.\\
Now $\neg\phi\times\psi\equiv\Box\bigl((\neg\phi|\psi)\leftrightarrow\neg\phi\bigr)\equiv\Box\bigl(\neg(\phi|\psi)\leftrightarrow\neg\phi\bigr)\equiv
\Box\bigl((\phi|\psi)\leftrightarrow\phi\bigr)\equiv\phi\times\psi$ by b4.\\
By applying b5 again, it comes $\psi\times\neg\phi\equiv\psi\times\phi$\,.
\paragraph{The empty universe.}\rien\\
It is deduced $\vdash\Box\neg\phi\rightarrow(\psi\times\neg\phi)$
and then $\vdash\Box\neg\phi\rightarrow(\psi\times\phi)$ by b5.weak.A.\\[5pt]
The remaining proof is obvious.
\paragraph{Left equivalences.}\rien\\[5pt]
\emph{Proof of the main theorem.}\\[5pt]
From $\vdash (\psi\rightarrow\eta)\rightarrow\bigl(\phi\rightarrow(\psi\rightarrow\eta)\bigr)$\,, it is deduced $\vdash\Box(\psi\rightarrow\eta)\rightarrow\Box\bigl(\phi\rightarrow(\psi\rightarrow\eta)\bigr)$\,.\\
By applying axiom b1, it comes $\vdash\Box(\psi\rightarrow\eta)\rightarrow\bigl(\Box\neg\phi\vee\Box(\psi\rightarrow\eta|\phi)\bigr)$\,.\\
Then axiom b2 implies $\vdash\Box(\psi\rightarrow\eta)\rightarrow\Bigl(\Box\neg\phi\vee\Box\bigl((\psi|\phi)\rightarrow(\eta|\phi)\bigr)\Bigr)$.\\
Since $\psi$ and $\eta$ are exchangeable, the theorem is deduced.
\\\\
\emph{Proof of the corollary.}\\[5pt]
It has been proved $\vdash\Box\neg\phi\rightarrow(\psi\times\phi)$, or equivalently $\vdash\Box\neg\phi\rightarrow\Box\bigl((\psi|\phi)\leftrightarrow\psi\bigr)$.\\
Of course, also holds $\vdash\Box\neg\phi\rightarrow\Box\bigl((\eta|\phi)\leftrightarrow\eta\bigr)$.\\
Then $\vdash\bigl(\Box\neg\phi\wedge\Box(\psi\leftrightarrow\eta)\bigl)\rightarrow\Box\bigl((\psi|\phi)\leftrightarrow(\eta|\phi)\bigr)$.\\
The corollary is then deduced from the main proposition.
\paragraph{Sub-universes are classical.}\rien\\
The first theorem is a consequence of axiom b4.
\vspace{5pt}\\
From axiom b2\,, it is deduced $\vdash(\psi\rightarrow\neg\eta|\phi)\rightarrow\bigl((\psi|\phi)\rightarrow(\neg\eta|\phi)\bigr)$\,.\\
It is deduced $\vdash\bigl((\psi|\phi)\wedge\neg(\neg\eta|\phi)\bigr)\rightarrow\neg(\psi\rightarrow\neg\eta|\phi)$\,.\\
Applying b4, it comes $\vdash\bigl((\psi|\phi)\wedge(\eta|\phi)\bigr)\rightarrow(\psi\wedge\eta|\phi)$.
\\[5pt]
Now $\vdash\phi\rightarrow\bigl((\psi\wedge\eta)\rightarrow\psi\bigr)$ and b1 imply $\vdash\Box\neg\phi\vee\Box\bigl((\psi\wedge\eta)\rightarrow\psi\big|\phi\bigr)$.\\
By b2 it is deduced $\vdash\Box\neg\phi\vee\Box\bigl((\psi\wedge\eta|\phi)\rightarrow(\psi|\phi)\bigr)$.\\
It is similarly proved $\vdash\Box\neg\phi\vee\Box\bigl((\psi\wedge\eta|\phi)\rightarrow(\eta|\phi)\bigr)$.\\
At last $\vdash\Box\neg\phi\vee\Box\Bigl((\psi\wedge\eta|\phi)\rightarrow\bigl((\psi|\phi)\wedge(\eta|\phi)\bigr)\Bigr)$.\\[3pt]
Now it has been shown $\vdash\Box\neg\phi\rightarrow\Box\bigl((\Xi|\phi)\leftrightarrow\Xi\bigr)$, and considering $\Xi=\psi$, $\eta$ or $\psi\wedge\eta$, it is implied $\vdash\Box\neg\phi\rightarrow\Box\Bigl((\psi\wedge\eta|\phi)\rightarrow\bigl((\psi|\phi)\wedge(\eta|\phi)\bigr)\Bigr)$.\\
At last $\vdash\Box\Bigl((\psi\wedge\eta|\phi)\rightarrow\bigl((\psi|\phi)\wedge(\eta|\phi)\bigr)\Bigr)$ and then $\vdash(\psi\wedge\eta|\phi)\rightarrow\bigl((\psi|\phi)\wedge(\eta|\phi)\bigr)$.\\
The second theorem is then proved.
~\vspace{5pt}\\
Third theorem is a consequence of the first and second theorems.
~\vspace{5pt}\\
Last theorem is a consequence of the first and third theorems.
\paragraph{Evaluating $(\top|\cdot)$ and $(\bot|\cdot)$\,.}\rien\\
From $\vdash\Box\psi\rightarrow\Box(\phi\rightarrow\psi)$ and b1, it comes $\vdash\Box\psi\rightarrow\bigl(\Box\neg\phi\vee\Box(\psi|\phi)\bigr)$\,.\\
Now $\vdash\Box\neg\phi\rightarrow\Box\bigl((\psi|\phi)\leftrightarrow\psi\bigr)$\,,
and consequently $\vdash\Box\psi\rightarrow\Box(\psi|\phi)$\,.
\paragraph{Inference property.}\rien\\
From b3 it comes $\vdash(\neg\psi|\phi)\rightarrow(\phi\rightarrow\neg\psi)$\,.\\
Then $\vdash\neg(\phi\rightarrow\neg\psi)\rightarrow\neg(\neg\psi|\phi)$ and 
$\vdash(\phi\wedge\psi)\rightarrow(\psi|\phi)$\,.\\
At last $\vdash(\phi\wedge\psi)\rightarrow\bigl((\psi|\phi)\wedge\phi\bigr)$\,.\\[5pt]
Conversely $\vdash(\psi|\phi)\rightarrow(\phi\rightarrow\psi)$ implies $\vdash\bigl((\psi|\phi)\wedge\phi\bigr)\rightarrow\bigl((\phi\rightarrow\psi)\wedge\phi\bigr)$\,.\\
Since $(\phi\rightarrow\psi)\wedge\phi\equiv\phi\wedge\psi$\,, the converse is proved.
\paragraph{Introspection.}\rien\\
Obvious from $\vdash \phi\rightarrow\phi$ and b1\,.
\paragraph{Inter-independence.}\rien\\
It is proved:
$$
\bigl((\psi|\phi)\big|\phi\bigr)\wedge(\phi|\phi)
\equiv\bigl((\psi|\phi)\wedge\phi\big|\phi\bigr)\equiv(\phi\wedge\psi|\phi)\equiv(\psi|\phi)\wedge(\phi|\phi)\;.
$$
As a consequence $\vdash(\phi|\phi)\rightarrow\Bigl(\bigl((\psi|\phi)\big|\phi\bigr)\leftrightarrow(\psi|\phi)\Bigr)$\,.\\
Then $\vdash\Box(\phi|\phi)\rightarrow\bigl((\psi|\phi)\times\phi\bigr)$\,.\\
Now $\vdash\Box\neg\phi\vee\Box(\phi|\phi)$ and $\vdash\Box\neg\phi\rightarrow\bigl((\psi|\phi)\times\phi\bigr)$ from previous results.\\
At last $\vdash(\psi|\phi)\times\phi$\,.
\paragraph{Independence invariance.}\rien\\
First theorem comes from the deduction:
$$
(\psi\times\phi)\equiv\Box\bigl((\psi|\phi)\leftrightarrow\psi\bigr)\equiv\Box\bigl(\neg(\psi|\phi)\leftrightarrow\neg\psi\bigr)\equiv\Box\bigl((\neg\psi|\phi)\leftrightarrow\neg\psi\bigr)\equiv(\neg\psi\times\phi)\;.
$$
The second theorem is also derived from similar deductions:
$$
\vdash\Box\Bigl(\bigl((\psi|\phi)\leftrightarrow\psi\bigr)\wedge\bigl((\eta|\phi)\leftrightarrow\eta\bigr)\Bigr)
\rightarrow
\Box\Bigl(\bigl((\psi|\phi)\wedge(\eta|\phi)\bigr)\leftrightarrow(\psi\wedge\eta)\Bigr)\vspace{-5pt}
$$
and then $\vdash\Box\Bigl(\bigl((\psi|\phi)\leftrightarrow\psi\bigr)\wedge\bigl((\eta|\phi)\leftrightarrow\eta\bigr)\Bigr)
\rightarrow
\Box\bigl((\psi\wedge\eta|\phi)\bigr)\leftrightarrow(\psi\wedge\eta)\bigr)$\,.\vspace{5pt}
\\
Now, let prove the third theorem.\\
The \emph{Left equivalences} theorem implies $\vdash\Box(\psi\leftrightarrow\eta)\rightarrow\Box\Bigl(\bigl((\psi|\phi)\leftrightarrow(\eta|\phi)\bigr)\wedge(\psi\leftrightarrow\eta)\Bigr)$.\\
Since $\vdash\bigl((\alpha\leftrightarrow\beta)\wedge(\gamma\leftrightarrow\delta)\bigr)\rightarrow\bigl((\alpha\leftrightarrow\gamma)\leftrightarrow(\beta\leftrightarrow\delta)\bigr)$, it is deduced\\ \rien\hspace{150pt}$\vdash\Box(\psi\leftrightarrow\eta)\rightarrow\Box\Bigl(\bigl((\psi|\phi)\leftrightarrow\psi\bigr)\leftrightarrow\bigl((\eta|\phi)\leftrightarrow\eta\bigr)\Bigr)$.\\
At last $\vdash\Box(\psi\leftrightarrow\eta)\rightarrow\Bigl(\Box\bigl((\psi|\phi)\leftrightarrow\psi\bigr)\leftrightarrow\Box\bigl((\eta|\phi)\leftrightarrow\eta\bigr)\Bigr)$.
\paragraph{Narcissistic independence.}\rien\\
From $\vdash\phi\rightarrow\phi$ it is deduced $\vdash\Box\neg\phi\vee\Box(\phi|\phi)$\,.\\
From definition, it is derived $\vdash(\phi\times\phi)\rightarrow\Box\bigl((\phi|\phi)\rightarrow\phi\bigr)$ and then $\vdash(\phi\times\phi)\rightarrow\bigl(\Box(\phi|\phi)\rightarrow\Box\phi\bigr)$\,.\\
It is thus deduced $\vdash(\phi\times\phi)\rightarrow(\Box\neg\phi\vee\Box\phi)$\,.
\paragraph{Independence and proof.}\rien\\
Combining $\vdash(\phi\vee\psi)\rightarrow(\neg\phi\rightarrow\psi)$ with b1 implies $\vdash\Box(\phi\vee\psi)\rightarrow\bigl(\Box\phi\vee\Box(\psi|\neg\phi)\bigr)$\,.\\
From b5.weak.A, it comes $\vdash(\psi\times\phi)\rightarrow\Box\bigl((\psi|\neg\phi)\leftrightarrow\psi\bigr)$\,.\\
As a consequence $\vdash(\phi\times\psi)\rightarrow\bigl(\Box(\phi\vee\psi)\rightarrow(\Box\phi\vee\Box\psi)\bigr)$\,.
\paragraph{Independence and regularity.}\rien\\[5pt]
\emph{Proof of the main theorem.}\\[5pt]
It is easy to prove $(\phi\wedge\eta)\rightarrow(\psi\wedge\eta)\equiv\neg\eta\vee(\phi\rightarrow\psi)$\,.\\
Then  $\vdash\Box\bigl((\phi\wedge\eta)\rightarrow(\psi\wedge\eta)\bigr)\rightarrow\Box\bigl(\neg\eta\vee(\phi\rightarrow\psi)\bigr)$\,.\\
Now $\vdash\bigl((\phi\times\eta)\wedge(\psi\times\eta)\bigr)\rightarrow\bigl((\phi\rightarrow\psi)\times\neg\eta\bigr)$\,,
by independence invariance and b5.weak.A.\\
The proof is achieved by means of the preceding property, \emph{independence and proof}.
\\\\
\emph{Proof of Corollary 2.}\\[5pt]
Assume $\vdash X\times\phi$ and $\vdash\Diamond\phi$.\\
Since $\vdash(\psi|\phi)\times\phi$ and $\psi\wedge\phi\equiv(\psi|\phi)\wedge\phi$, it is deduced from $X\wedge\phi\equiv\psi\wedge\phi$ that $X\equiv(\psi|\phi)$.
\paragraph{Right equivalences.}\rien\\
First notice that all previous properties are obtained without b5.weak.B.
\\[5pt]
From $(\phi|\psi)\wedge\psi\equiv\phi\wedge\psi$, $(\phi|\eta)\wedge\eta\equiv\phi\wedge\eta$ and $\vdash(\psi\leftrightarrow\eta)\rightarrow\bigl((\phi\wedge\psi)\leftrightarrow(\phi\wedge\eta)\bigr)$, it is deduced $\vdash(\psi\leftrightarrow\eta)\rightarrow\Bigl(\bigl((\phi|\psi)\wedge\psi\bigr)\leftrightarrow\bigl((\phi|\eta)\wedge\eta\bigr)\Bigr)\,.$\\
Then $\vdash(\psi\leftrightarrow\eta)\rightarrow\Bigl(\bigl((\phi|\psi)\wedge\psi\bigr)\leftrightarrow\bigl((\phi|\eta)\wedge\psi\bigr)\Bigr)$ and finally:\vspace{-5pt}
$$
\vdash\Box(\psi\leftrightarrow\eta)\rightarrow\Box\Bigl(\bigl((\phi|\psi)\wedge\psi\bigr)\leftrightarrow\bigl((\phi|\eta)\wedge\psi\bigr)\Bigr)\,.
$$
Now $\vdash(\phi|\psi)\times\psi$ and $\vdash(\phi|\eta)\times\eta$\,.
\\
Since $\vdash\Box(\psi\leftrightarrow\eta)\rightarrow\Bigr(\bigl(\eta\times(\phi|\eta)\bigr)\leftrightarrow\bigl(\psi\times(\phi|\eta)\bigr)\Bigl)$, it comes $\vdash\Box(\psi\leftrightarrow\eta)\rightarrow\bigl((\phi|\eta)\times\psi\bigr)$ by b5.
\\
Finally $\vdash\Box(\psi\leftrightarrow\eta)\rightarrow\biggl(\bigl((\phi|\psi)\times\psi\bigr)\wedge\bigl((\phi|\eta)\times\psi\bigr)\wedge\Box\Bigl(\bigl((\phi|\psi)\wedge\psi\bigr)\leftrightarrow\bigl((\phi|\eta)\wedge\psi\bigr)\Bigr)\biggr)$\,.\\
Applying the regularity, it comes $\vdash\Box(\psi\leftrightarrow\eta)\rightarrow\Bigl(\Box\neg\psi\vee\Box\bigl((\phi|\psi)\leftrightarrow(\phi|\eta)\bigr)\Bigr)$\,.\\[5pt]
Now $\vdash\Box(\psi\leftrightarrow\eta)\rightarrow (\Box\neg\psi\leftrightarrow\Box\neg\eta)$
and $\vdash\Box\neg\Xi\rightarrow\Box\bigl((\phi|\Xi)\leftrightarrow\phi\bigr)$ for $\Xi=\psi$ or $\eta$.\\
It is deduced $\vdash\bigl(\Box\neg\psi\wedge\Box(\psi\leftrightarrow\eta)\bigr)\rightarrow\Box\bigl((\phi|\psi)\leftrightarrow(\phi|\eta)\bigr)$\,, thus completing the proof.
\paragraph{Reduction rule.}\rien\\
Since $\vdash(\psi|\phi)\times\phi$, b5 implies $\vdash\phi\times(\psi|\phi)$ and $\vdash \Box\Bigl(\bigl(\phi\big|(\psi|\phi)\bigr)\leftrightarrow\phi\Bigr)$\,.
\paragraph{Markov Property.}\rien\\
Since $\vdash(\phi_t|\phi_{t-1})\times\phi_{t-1}$, it comes $\vdash\left(
\bigwedge_{\tau=1}^{t-2}\bigl((\phi_t|\phi_{t-1})\times\phi_{\tau}\bigr)
\right)\rightarrow\left(
\bigwedge_{\tau=1}^{t-1}\bigl((\phi_t|\phi_{t-1})\times\phi_{\tau}\bigr)
\right)$\,.\\
Then $\vdash\left(
\bigwedge_{\tau=1}^{t-2}\bigl((\phi_t|\phi_{t-1})\times\phi_{\tau}\bigr)
\right)\rightarrow\left((\phi_t|\phi_{t-1})\times \left(\bigwedge_{\tau=1}^{t-1}\phi_\tau\right)\right)$\,.\\
Now, $(\phi_t|\phi_{t-1})\wedge \left(\bigwedge_{\tau=1}^{t-1}\phi_\tau\right)\equiv\bigwedge_{\tau=1}^{t}\phi_\tau
\equiv\left(\phi_t\left|\bigwedge_{\tau=1}^{t-1}\phi_\tau\right.\right)\wedge \left(\bigwedge_{\tau=1}^{t-1}\phi_\tau\right)$\,.\\
Since $\vdash\left(\phi_t\left|\bigwedge_{\tau=1}^{t-1}\phi_\tau\right.\right)\times \left(\bigwedge_{\tau=1}^{t-1}\phi_\tau\right)$\,, the proof is achieved by applying the regularity.
\paragraph{Link between $\bigl((\eta|\psi)\big|\phi\bigr)$ and $(\eta|\phi\wedge\psi)$.}\rien\\[5pt]
\emph{Proof of the logical counterpart to Lewis' triviality.}
\\[5pt]
Since $\neg(\phi\rightarrow\psi)\equiv\phi\wedge\neg\psi$\,, it is equivalent to prove:
$$
\vdash\Bigl(\Diamond(\phi\wedge\psi)\wedge\bigl(\Diamond(\neg\psi\wedge\phi)\vee\Diamond(\neg\phi\wedge\psi)\bigr)\Bigr)\rightarrow(\phi\times\psi)
\;.
$$
Since $\times$ is symmetric, it is sufficient to prove $\vdash\bigl(\Diamond(\phi\wedge\psi)\wedge\Diamond(\neg\psi\wedge\phi)\bigr)\rightarrow(\phi\times\psi)$\,.\\[5pt]
The \emph{introspection} property implies $\vdash\Diamond(\phi\wedge\psi)\rightarrow\Box(\phi\wedge\psi|\phi\wedge\psi)$, denoted $(a)$, and $\vdash\Diamond(\neg\psi\wedge\phi)\rightarrow\Box(\neg\psi\wedge\phi|\neg\psi\wedge\phi)$\,. \\
It is thus deduced $\vdash\Diamond(\phi\wedge\psi)\rightarrow\Box\Bigr((\psi|\psi\wedge\phi)\leftrightarrow\bigr((\psi|\psi\wedge\phi)\wedge(\psi\wedge\phi|\psi\wedge\phi)\bigl)\Bigl)$, denoted $(b)$, and $\vdash\Diamond(\neg\psi\wedge\phi)\rightarrow\Box\Bigr((\psi|\neg\psi\wedge\phi)\leftrightarrow\bigr((\psi|\neg\psi\wedge\phi)\wedge(\neg\psi\wedge\phi|\neg\psi\wedge\phi)\bigl)\Bigl)$. \\
From the deduction $(\psi|\neg\psi\wedge\phi)\wedge(\neg\psi\wedge\phi|\neg\psi\wedge\phi)\equiv(\bot|\neg\psi\wedge\phi)\equiv\bot$, it is derived $\vdash\Diamond(\neg\psi\wedge\phi)\rightarrow\Box\bigl((\psi|\neg\psi\wedge\phi)\leftrightarrow\bot\bigr)$, denoted $(c)$. \\
From the deduction $(\psi|\psi\wedge\phi)\wedge(\psi\wedge\phi|\psi\wedge\phi)\equiv(\psi\wedge\phi|\psi\wedge\phi)$, $(a)$ and $(b)$, it comes $\vdash\Diamond(\phi\wedge\psi)\rightarrow\Box\bigl((\psi|\psi\wedge\phi))\leftrightarrow\top\bigr)$, denoted $(d)$.\\
Now $(\psi|\phi)\equiv\bigl((\psi|\phi)\wedge\psi\bigr)\vee\bigl((\psi|\phi)\wedge\neg\psi\bigr)\equiv
\Bigl(\bigl((\psi|\phi)\big|\psi\bigr)\wedge\psi\Bigr)\vee\Bigl(\bigl((\psi|\phi)\big|\neg\psi\bigr)\wedge\neg\psi\Bigr)$\,, and by applying axiom $(\ast)$, $(\psi|\phi)\equiv\bigl((\psi|\phi\wedge\psi)\wedge\psi\bigr)\vee\bigl((\psi|\phi\wedge\neg\psi)\wedge\neg\psi\bigr)$.\\
Then $\vdash\bigl(\Diamond(\neg\psi\wedge\phi)\wedge\Diamond(\phi\wedge\psi)\bigr)\rightarrow\Box\Bigl(
(\psi|\phi)\leftrightarrow\bigl((\top\wedge\psi)\vee(\bot\wedge\neg\psi)\bigr)
\Bigr)$ by $(c)$ and $(d)$. \\
At last $\vdash\bigl(\Diamond(\neg\psi\wedge\phi)\wedge\Diamond(\phi\wedge\psi)\bigr)\rightarrow\Box\bigl(
(\psi|\phi)\leftrightarrow\psi\bigr)$\,.
\section{Proof: model transfer}
\label{proof:modtrans}
First notice that the above construction of $H$ is possible for any proposition $\phi\in\mathcal{L}$, since it is always obtained $H(\phi)\in M$\,.\\
Now, $R=W\times W$ implies $H(\Box\phi)=\bigl\{t\in W\big/ \forall u\in W,\,(t,u)\in R\Rightarrow u\in H(\phi)\bigr\}$\,, so that $(W,R,H,f)$ is actually a Kripke model.\\
Let verify the compliance with m3, b1, b2, b3, b4 and b5 (resp. b5.weak.$\ast$).\\[10pt]$\bullet$
\underline{Compliance with m3} is obtained from the fact that $R$ is reflexive.
\\[5pt]$\bullet$
Proof of $H\Bigl(\Box(\phi\rightarrow\psi)\rightarrow\bigl(\Box\neg\phi\vee\Box(\psi|\phi)\bigr)\Bigr)=W$\,, \emph{i.e.} \underline{compliance with b1}.\\
By definition, $H\Bigl(\Box(\phi\rightarrow\psi)\rightarrow\bigl(\Box\neg\phi\vee\Box(\psi|\phi)\bigr)\Bigr)=H\bigl(\Box(\psi|\phi)\bigr)\cup H(\Box\neg\phi) \cup \Bigl(W\setminus H\bigl(\Box(\phi\rightarrow\psi)\bigr)\Bigl)$\,.\\
Cases $H\bigl(\Box(\phi\rightarrow\psi)\bigr)=\emptyset$ or $H(\Box\neg\phi)=W$ then imply \mbox{$H\Bigl(\Box(\phi\rightarrow\psi)\rightarrow\bigl(\Box\neg\phi\vee\Box(\psi|\phi)\bigr)\Bigr)=W$\,.\rien}\\
Otherwise $H\bigl(\Box(\phi\rightarrow\psi)\bigr)\ne\emptyset$ and $H(\Box\neg\phi)\ne W$\,, thus implying $H(\phi\rightarrow\psi)=W$ and $H(\neg\phi)\ne W$\,.\\
It is deduced $H(\phi)\subset H(\psi)$ and $H(\phi)\ne\emptyset$\,, and then $f\bigl(H(\psi), H(\phi)\bigr)=W$ by using $\beta1$.\\
Finally $H\bigl(\Box(\psi|\phi)\bigr)=W$ and again $H\Bigl(\Box(\phi\rightarrow\psi)\rightarrow\bigl(\Box\neg\phi\vee\Box(\psi|\phi)\bigr)\Bigr)=W$\,.
\\[5pt]$\bullet$
Proof of $H\bigl(\neg(\neg\psi|\phi)\leftrightarrow(\psi|\phi)\bigr)=W$\,, \emph{i.e.} \underline{compliance with b4}.\\
It is deduced $
H\bigl((\neg\psi|\phi)\bigr)=f\bigl(H(\neg\psi),H(\phi)\bigr)=f\bigl(W\setminus H(\psi),H(\phi)\bigr)=W\setminus f\bigl(H(\psi),H(\phi)\bigr)=W\setminus H\bigl((\psi|\phi)\bigr)
$\,, by using $\beta4$.\\
Then $\Bigl(H\bigl((\neg\psi|\phi)\bigr)\cap \bigl(W\setminus H\bigl((\psi|\phi)\bigr)\bigr)\Bigr)\cup\Bigl(\bigl(W\setminus H\bigl((\neg\psi|\phi)\bigr)\bigr)\cap H\bigl((\psi|\phi)\bigr)\Bigr)=W$\,.\\
And finally $H\bigl(\neg(\neg\psi|\phi)\leftrightarrow(\psi|\phi)\bigr)=W$\,.
\\[5pt]$\bullet$
Proof of $H\Bigl((\psi\rightarrow\eta|\phi)\rightarrow\bigl((\psi|\phi)\rightarrow(\eta|\phi)\bigr)\Bigr)=W$\,, \emph{i.e.} \underline{compliance with b2}.\\
By $\beta2$, it is deduced
$H\bigl((\psi\rightarrow\eta|\phi)\bigr)=f\bigl(H(\psi\rightarrow\eta),H(\phi)\bigr)=f\bigl(H(\neg\psi)\cup H(\eta),H(\phi)\bigr)$\\
\rien\hspace{25pt}$\subset f\bigl(H(\neg\psi),H(\phi)\bigr)\cup f\bigl(H(\eta),H(\phi)\bigr)=H\bigl((\neg\psi|\phi)\bigr)\cup H\bigl((\eta|\phi)\bigr) =H\bigl(\neg(\psi|\phi)\bigr)\cup H\bigl((\eta|\phi)\bigr)$\,.\\
Then $H\Bigl((\psi\rightarrow\eta|\phi)\rightarrow\bigl((\psi|\phi)\rightarrow(\eta|\phi)\bigr)\Bigr)= W$\,.
\\[5pt]$\bullet$
Proof of $H\bigl((\psi|\phi)\rightarrow(\phi\rightarrow\psi)\bigr)=W$\,, \emph{i.e.} \underline{compliance with b3}.\\
Immediate from $\beta3$, \emph{i.e.} $H(\phi)\cap f\bigl(H(\psi),H(\phi)\bigr)\subset H(\psi)$\,.
\\[5pt]
{\bf Case DmBL.}\\$\bullet$
Proof of $H\Bigl(\Box\bigl((\psi|\phi)\leftrightarrow\psi\bigr)\leftrightarrow\Box\bigl((\phi|\psi)\leftrightarrow\phi\bigr)\Bigr)=W$\,, \emph{i.e.} \underline{compliance with b5}.\\
By $\beta5$, $f\bigl(H(\psi),H(\phi)\bigr)=H(\psi)$ if and only if $f\bigl(H(\phi),H(\psi)\bigr)=H(\phi)$\,.\\
Then $H\bigl((\psi|\phi)\leftrightarrow\psi\bigr)=W$ if and only if $H\bigl((\phi|\psi)\leftrightarrow\phi\bigr)=W$\,.\\
As a consequence, $H\Bigl(\Box\bigl((\psi|\phi)\leftrightarrow\psi\bigr)\Bigr)=H\Bigl(\Box\bigl((\phi|\psi)\leftrightarrow\phi\bigr)\Bigr)$ and the result.
\\[5pt]
{\bf Case DmBL$_\ast$.}\\$\bullet$
Proof of $H\Bigl(\Box\bigl((\psi|\phi)\leftrightarrow\psi\bigr)\leftrightarrow\Box\bigl((\psi|\neg\phi)\leftrightarrow\psi\bigr)\Bigr)=W$\,, \emph{i.e.} \underline{compliance with b5.weak.A}.\\
By $\beta5w$, $f\bigl(H(\psi),H(\phi)\bigr)=H(\psi)$ if and only if $f\bigl(H(\psi),W\setminus H(\phi)\bigr)=H(\psi)$\,.\\
Then $H\bigl((\psi|\phi)\leftrightarrow\psi\bigr)=W$ if and only if $H\bigl((\psi|\neg\phi)\leftrightarrow\psi\bigr)=W$\,.\\
As a consequence, $H\Bigl(\Box\bigl((\psi|\phi)\leftrightarrow\psi\bigr)\Bigr)=H\Bigl(\Box\bigl((\psi|\neg\phi)\leftrightarrow\psi\bigr)\Bigr)$ and the result.
\\[5pt]$\bullet$
Proof of $H\Bigl(\Box(\psi\leftrightarrow\eta)\rightarrow\Box\bigl((\phi|\psi)\leftrightarrow(\phi|\eta)\bigr)\Bigr)=W$\,, \emph{i.e.} \underline{compliance with b5.weak.B}.\\
Assume first $H(\psi)\ne H(\eta)$\,.\\
It comes $H(\psi\leftrightarrow\eta)\ne W$ and $H\bigl(\Box(\psi\leftrightarrow\eta)\bigr)=\emptyset$\,.\\[3pt]
Assume now $H(\psi)= H(\eta)$\,.\\
Then $H\bigl((\phi|\psi)\bigr)=H\bigl((\phi|\eta)\bigr)$ and $H\Bigl(\Box\bigl((\phi|\psi)\leftrightarrow(\phi|\eta)\bigr)\Bigr)=H\bigl((\phi|\psi)\leftrightarrow(\phi|\eta)\bigr)=W$\,.\\[3pt]
At last, both cases imply the result.
%
%
\section{Proof: properties of $(\Omega_n, M_n, h_n, f_n,\Lambda_n,\mu_n)_{n\in\Nset}$}
\label{Appendix:MainProof}
To be proved:
\emph{\begin{description}
\item[$\rien\quad\bullet_\mu$] $\mu_n:M_n\rightarrow M_{n+1}$ is a one-to-one Boolean morphism\,,
\item[$\rien\quad\bullet_f$] $f_{n+1}\bigl(\mu_n(B),\mu_n(A)\bigr)=\mu_n\bigl(f_n(B,A)\bigr)$\,,
\item[$\rien\quad\tilde\beta1$.] $A\subset B$ and $A\ne\emptyset$ imply $f_n(B,A)=\Omega_n$\,,
\item[$\rien\quad\tilde\beta2$.] $f_n(B\cup C,A)= f_n(B,A)\cup f_n(C,A)$\,,
\item[$\rien\quad\tilde\beta3$.] $A\cap f_n(B,A) = A\cap B$\,,
\item[$\rien\quad\tilde\beta4$.] $f_n(\sim B,A)=\sim f_n(B,A)$\,,
\item[$\rien\quad\tilde\beta5w$.] $f_n(B,A)=B$ implies $f_n(B,\sim A)=B$\,,
\end{description}
being assumed $A,B,C\in M_n$\,, and $f_n(\cdot,\cdot)$ defined for the considered cases.
}
\\[5pt]
The proof is recursive and needs to consider the two cases in the definition of $(\mu_n,f_n)$\,.
\\[3pt]
The properties $\tilde\beta\ast$ are obvious for $n=0$, since $f_0$ is only defined by $f_(A,\emptyset)=f_0(A,\Omega_0)=A$\,.
From now on, it is assumed that $\tilde\beta\ast$ hold true for $k\le n$, and that $\bullet_\mu$ and $\bullet_f$ hold true for $k\le n-1$\,.
The subsequent paragraphs establish the proof of $\beta\ast$ for $n+1$ and the proof of $\bullet_\mu$ and $\bullet_f$ for $n$.
\paragraph{Preliminary remark.} It is noticed that $\tilde\beta2$ and $\tilde\beta4$ imply:
$$\tilde\beta6:\ f_n(B\cap C,A)= f_n(B,A)\cap f_n(C,A)\,.$$
\subsection{Lemma.}
\label{lemma:include}
$\bigcup_{i\in I_n}\Pi_n(i)=b_n$
and
$\bigcup_{i\in I_n}\Gamma_n(i)=\sim b_n$\,;
in particular, $\Pi_n(i)\cap\Gamma_n(j)=\emptyset$ for any $i,j\in I_n$\,.\\[5pt]
Moreover $\Pi_n(i)\cap\Pi_n(j)=\Gamma_n(i)\cap\Gamma_n(j)=\emptyset$ for any $i,j\in I_n$ such that $i\ne j$\,.
\begin{description}
\item[Proof.]The proof is obvious for case 1.\\
Now, let consider case 0.\\
By definition $
\bigcup_{i\in I_n}\Pi_n(i)=\left(\bigcup_{\omega\in\mu_\nu(b_\nu)}\omega_{[n]}\right)
\cap
\left(\bigcup_{\omega'\in\sim\mu_\nu(b_\nu)}f_n(\omega'_{[n]},\sim b_n)\right)
$\,.\\
By recursion hypothesis over $\tilde\beta1$ it comes $f_n(\sim b_n,\sim b_n)=\Omega_n$\,.
\\
Then by $\tilde\beta2$\,, $\bigcup_{i\in I_n}\Pi_n(i)=b_n\cap f_n(\sim b_n,\sim b_n)=b_n\cap\Omega_n= b_n\;.$\\[5pt]
For any $\omega_1,\omega_2\in \sim\mu_\nu(b_\nu)$ such that $\omega_1\ne\omega_2$\,, it comes by $\tilde\beta6$ (deduced from $\tilde\beta2$ and $\tilde\beta4$)\,:
$$
f_n(\omega_{1[n]},\sim b_n)\cap f_n(\omega_{2[n]},\sim b_n)=f_n(\omega_{1[n]}\cap\omega_{2[n]},\sim b_n)=f_n(\emptyset,\sim b_n)=\emptyset\;.
$$
Finally $\Pi_n(i)\cap\Pi_n(j)=\emptyset$ for any $i,j\in I_n$ such that $i\ne j$\,.\\[5pt]
The results are similarly proved for $\Gamma_n$\,.
\item[$\Box\Box\Box$]\rien
\end{description}
\emph{Corollary 1.}
$$
\mu_n(b_n)=T(\sim\mu_n(b_n))=\bigcup_{i\in I_n}\Pi_n(i)\times\Gamma_n(i)
\quad\mbox{and}\quad
\sim\mu_n(b_n)=T(\mu_n(b_n))=\bigcup_{i\in I_n}\Gamma_n(i)\times\Pi_n(i)\;.
$$
\emph{Corollary 2.}
$$f_{n+1}\bigl(C,\mu_n(b_n)\bigr)=\bigl(C\cap\mu_n(b_n)\bigr)\cup\bigl(T(C)\cap\sim\mu_n(b_n)\bigr)$$
and
$$f_{n+1}\bigl(C,\sim\mu_n(b_n)\bigr)=\bigl(T(C)\cap\mu_n(b_n)\bigr)\cup\bigl(C\cap\sim\mu_n(b_n)\bigr)\;.$$
Both corollary are obvious from the definition.
\subsection{Proof of $\bullet_\mu$}
\label{Proof:Cap}
The following properties (whose proofs are immediate) will be useful:
\begin{description}
\item[$\ell1.$] $(A\cup B)\times C=(A\times C)\cup (B\times C)$
and $A\times (B\cup C)=(A\times B)\cup (A\times C)$\,, for any $A,B,C$\,,
\item[$\ell2.$] $(A\cap B)\times C=(A\times C)\cap (B\times C)$
and $A\times (B\cap C)=(A\times B)\cap (A\times C)$\,, for any $A,B,C$\,,
\item[$\ell3.$] $C\cap D=\emptyset$ implies $(C\times A)\cap(D\times B)=(A\times C)\cap(B\times D)=\emptyset$\,, for any $A,B,C,D$\,,
\item[$\ell4.$] $(A\cup B)\cap (C\cup D)=\emptyset$ implies $(A\cap B)\cup(C\cap D)=(A\cup C)\cap(B\cup D)$\,, for any $A,B,C,D$\,.
\item[$\ell5.$] $(A\cup B)\cap (C\cup D)=\emptyset$ and $A\cup C=B\cup D$ imply $A=B$ and $C=D$\,, for any $A,B,C,D$\,.
\item[$\ell6.$] $C\cap D=\emptyset$ implies $(C\cup D)\setminus\bigl((A\cap C)\cup(B\cap D)\bigr)=(C\setminus A)\cup(D\setminus B)$\,, for any $A,B,C,D$\,.
\end{description}
\paragraph{\underline{Proof of $\mu_n(\Omega_n)=\Omega_{n+1}$ and $\mu_n(\emptyset)=\emptyset$\,.}}
Immediate from the definitions.
\paragraph{\underline{Proof of $\mu_n(A\cap B)=\mu_n(A)\cap\mu_n(B)$\,.}}
By applying $\ell2$, it is deduced:
$$\begin{array}{@{}l@{}}
\mu_n(A\cap B)=\bigcup_{i\in I_n} \biggl(\Bigl(\bigl(A\cap\Pi_n(i)\bigr)\times\Gamma_n(i)\Bigr)
\cap
\Bigl(\bigl(B\cap\Pi_n(i)\bigr)\times\Gamma_n(i)\Bigr)\biggr)
\vspace{5pt}\\\rien\hspace{50pt}
\cup\ 
\bigcup_{i\in I_n}
\biggl(\Bigl(\bigl(A\cap\Gamma_n(i)\bigr)\times\Pi_n(i)\Bigr)
\cap
\Bigl(\bigl(B\cap\Gamma_n(i)\bigr)\times\Pi_n(i)\Bigr)\biggr)
\;.
\end{array}$$
By lemma~\ref{lemma:include}, and applying $\ell3$ and $\ell4$, it is deduced:
$$\begin{array}{@{}l@{}}
\mu_n(A\cap B)=\bigcup_{i\in I_n} \biggl(\Bigl(\bigl(A\cap\Pi_n(i)\bigr)\times\Gamma_n(i)\Bigr)
\cup
\Bigl(\bigl(A\cap\Gamma_n(i)\bigr)\times\Pi_n(i)\Bigr)\biggr)
\vspace{5pt}\\\rien\hspace{50pt}
\cap\ 
\bigcup_{i\in I_n}
\biggl(\Bigl(\bigl(B\cap\Pi_n(i)\bigr)\times\Gamma_n(i)\Bigr)
\cup
\Bigl(\bigl(B\cap\Gamma_n(i)\bigr)\times\Pi_n(i)\Bigr)\biggr)
=\mu_n(A)\cap\mu_n(B)\;.
\end{array}$$
\paragraph{\underline{Proof of $\mu_n(A\cup B)=\mu_n(A)\cup\mu_n(B)$\,.}}
Obviously deduced from $\ell1$.
\paragraph{\underline{$\mu_n$ is one-to-one.}}
Assume $\mu_n(A)=\mu_n(B)$\,;
then:
$$\begin{array}{@{}l@{}}
\bigcup_{i\in I_n} \biggl(\Bigl(\bigl(A\cap\Pi_n(i)\bigr)\times\Gamma_n(i)\Bigr)
\cup
\Bigl(\bigl(A\cap\Gamma_n(i)\bigr)\times\Pi_n(i)\Bigr)\biggr)
\vspace{5pt}\\\rien\hspace{50pt}
=
\bigcup_{i\in I_n}
\biggl(\Bigl(\bigl(B\cap\Pi_n(i)\bigr)\times\Gamma_n(i)\Bigr)
\cup
\Bigl(\bigl(B\cap\Gamma_n(i)\bigr)\times\Pi_n(i)\Bigr)\biggr)
\;.
\end{array}$$
By lemma~\ref{lemma:include}, and applying $\ell2$, $\ell3$ and $\ell5$, it is deduced for any $i\in I_n$\,:
$$
\bigl(A\cap\Pi_n(i)\bigr)\times\Gamma_n(i)
=
\bigl(B\cap\Pi_n(i)\bigr)\times\Gamma_n(i)
\ 
\mbox{ and }
\ 
\bigl(A\cap\Gamma_n(i)\bigr)\times\Pi_n(i)
=
\bigl(B\cap\Gamma_n(i)\bigr)\times\Pi_n(i)
\;.$$
Finally $A\cap\Pi_n(i)=B\cap\Pi_n(i)$ and $A\cap\Gamma_n(i)=B\cap\Gamma_n(i)$ for any $i\in I_n$\,, and:
$$
A\cap\bigcup_{i\in I_n}\bigl(\Pi_n(i)\cup\Gamma_n(i)\bigr)=B\cap\bigcup_{i\in I_n}\bigl(\Pi_n(i)\cup\Gamma_n(i)\bigr)\;.
$$
$A=B$ is deduced by applying the lemma.
\paragraph{Conclusion.} The previous results imply that $\mu_n$ is a one-to-one Boolean morphism.
\subsection{Proof of $\bullet_f$}\label{Proof:idempot}
By definition, the result holds true for any $A\in M_n\setminus\{\emptyset,\Omega_n,b_n,\sim b_n\}$\,.
It is also true for $A=\emptyset$ or $A=\Omega_n$\,, since $f_{n+1}\bigl(\mu_n(B),\mu_n(\emptyset)\bigl)=f_{n+1}\bigl(\mu_n(B),\emptyset\bigl)=\mu_n(B)=\mu_n\bigl(f_n(B,\emptyset)\bigr)$
and similarly
$f_{n+1}\bigl(\mu_n(B),\mu_n(\Omega_n)\bigl)=f_{n+1}\bigl(\mu_n(B),\Omega_{n+1}\bigl)=\mu_n(B)=\mu_n\bigl(f_n(B,\Omega_n)\bigr)$\,.
\\[5pt]
The true difficulties come from the cases $A=b_n$ or $A=\sim b_n$\,.\\
Subsequently, it is assumed $A=b_n$\,; the case $A=\sim b_n$ is quite similar.\\[5pt]
It comes:
$$
f_{n+1}\bigl(\mu_n(B),\mu_n(b_n)\bigr)=(\mathrm{id}\cup T)
\biggl(\mu_n(B)\cap\Bigl(\bigcup_{i\in I_n}\bigl(\Pi_n(i)\times\Gamma_n(i)\bigr)\Bigr)\biggr)
=(\mathrm{id}\cup T)\Bigl(\bigcup_{i\in I_n}\bigl(B\cap\Pi_n(i)\bigr)\times\Gamma_n(i)\Bigr)\;.
$$
The existence of $f_n(B,b_n)$ necessary implies the case 0\,,
and there is $C\in M_{\nu+1}$ such that $B=C_{[n]}$\,.\\
By recursion hypotheses $\bullet_\mu$, it comes $B\cap\omega_{[n]}=(C\cap\omega)_{[n]}=\omega_{[n]}$ if $\omega\in C$\,, $=\emptyset$ if $\omega\not\in C$\,.\\
Moreover, $f_n(\omega_{[n]},b_n)\cap f_n(B,b_n)=f_n(\omega_{[n]}\cap B,b_n)=f_n(\omega_{[n]},b_n)$ if $\omega\in C$\,, $=\emptyset$ if $\omega\not\in C$\,.\\
As a consequence $\bigcup_{i\in I_n}\bigl(B\cap\Pi_n(i)\bigr)\times\Gamma_n(i)=\bigcup_{i\in I_n}\Pi_n(i)\times\bigl(f_n(B,b_n)\cap\Gamma_n(i)\bigr)$\,.
\\[3pt]
By $\tilde\beta3$,
$B\cap\omega_{[n]}=B\cap b_n\cap\omega_{[n]}=f_n(B,b_n)\cap b_n\cap\omega_{[n]}=f_n(B,b_n)\cap\omega_{[n]}$ for any $\omega\in\mu_\nu(b_\nu)$\,.\\
As a consequence $\bigcup_{i\in I_n}\bigl(B\cap\Pi_n(i)\bigr)\times\Gamma_n(i)=\bigcup_{i\in I_n}\bigl(f_n(B,b_n)\cap\Pi_n(i)\bigr)\times\Gamma_n(i)$\,.
\\
By applying the both results, it comes:
$$
f_{n+1}\bigl(\mu_n(B),\mu_n(b_n)\bigr)=
\Bigl(\bigcup_{i\in I_n}\bigl(f_n(B,b_n)\cap\Pi_n(i)\bigr)\times\Gamma_n(i)\Bigr)
\cup
\Bigl(\bigcup_{i\in I_n}\bigl(f_n(B,b_n)\cap\Gamma_n(i)\bigr)\times\Pi_n(i)\Bigr)\;.
$$
And by definition of $\mu_n$, it is finally deduced $f_{n+1}\bigl(\mu_n(B),\mu_n(b_n)\bigr)=\mu_n\bigl(f_{n}(B,b_n)\bigr)$\,.
\subsection{Proof of $\tilde\beta1$}
For $A\not\in \{\mu_n(b_n),\sim \mu_n(b_n),\emptyset,\Omega_{n+1}\}$, the propriety is inherited from $n$ by applying $\bullet_f$\,.\\
The property is also obvious for $A\in \{\emptyset,\Omega_{n+1}\}$\,.
\\
The difficulty comes from $A=\mu_n(b_n)$ or $A=\sim \mu_n(b_n)$\,;
then notice that $A\ne\emptyset$ by construction.
\\[5pt]
It is now hypothesized $A=\mu_n(b_n)\subset B$\,;
the case $A=\sim \mu_n(b_n)$ is quite similar.\\
Then $T\bigl(\mu_n(b_n)\bigr)\subset T(B)$ and by lemma, corollary 1\&2:
$$
f_{n+1}\bigl(B,\mu_n(b_n)\bigr)=
\bigl(B\cap\mu_n(b_n)\bigr)\cup\bigl(T(B)\cap\sim\mu_n(b_n)\bigr)=\mu_n(b_n)\cup\sim\mu_n(b_n)
=\Omega_{n+1}\;.
$$
\subsection{Proof of $\tilde\beta2$}
For $A\not\in \{\mu_n(b_n),\sim \mu_n(b_n),\emptyset,\Omega_{n+1}\}$, the propriety is inherited from $n$ by applying $\bullet_f$\,.\\
The property is also obvious for $A\in \{\emptyset,\Omega_{n+1}\}$\,.
\\
The property is then immediate for $A\in\{\mu_n(b_n),\sim \mu_n(b_n)\}$\,, since $T(B_1\cup B_2)=T(B_1)\cup T(B_2)$.
\subsection{Proof of $\tilde\beta3$}
For $A\not\in \{\mu_n(b_n),\sim \mu_n(b_n),\emptyset,\Omega_{n+1}\}$, the propriety is inherited from $n$ by applying $\bullet_f$\,.
The property is also obvious for $A\in \{\emptyset,\Omega_{n+1}\}$\,.
\\
The difficulty comes from $A\in \{\mu_n(b_n),\sim \mu_n(b_n)\}$.
\\[5pt]
It is now hypothesized $A=\mu_n(b_n)$\,;
the case $A=\sim \mu_n(b_n)$ is quite similar.\\
The result is immediate from corollary 2 of lemma.
\subsection{Proof of $\tilde\beta4$}
For $A\not\in \{\mu_n(b_n),\sim \mu_n(b_n),\emptyset,\Omega_{n+1}\}$, the propriety is inherited from $n$ by applying $\bullet_f$\,.
The property is also obvious for $A\in \{\emptyset,\Omega_{n+1}\}$\,.
\\
The difficulty comes from $A\in \{\mu_n(b_n),\sim \mu_n(b_n)\}$\,.\\[5pt]
It is now hypothesized $A=\mu_n(b_n)$\,;
the case $A=\sim \mu_n(b_n)$ is quite similar.\\
By corollary 2 of lemma, $f_{n+1}\bigl(\sim B,\mu_{n}(b_n)\bigr)=\bigl(\mu_n(b_n)\setminus B\bigr)\cup\bigl(\sim\mu_n(b_n)\setminus T(B)\bigr)$\,.\\
By $\ell6$, $f_{n+1}\bigl(\sim B,\mu_{n}(b_n)\bigr)=
\sim\Bigl(\bigl(B\cap\mu_n(b_n)\bigr)\cup\bigl(T(B)\cap\sim\mu_n(b_n)\bigr)\Bigr)=\sim f_{n+1}\bigl(B,\mu_{n}(b_n)\bigr)$\,.
\subsection{Lemma 2.}
Let $C\in M_{n+1}$\,. Then:
$$
f_{n+1}\Bigl(f_{n+1}\bigl(C,\mu_n(b_n)\bigr),\mu_n(b_n)\Bigr)=
f_{n+1}\Bigl(f_{n+1}\bigl(C,\mu_n(b_n)\bigr),\sim\mu_n(b_n)\Bigr)=
f_{n+1}\bigl(C,\mu_n(b_n)\bigr)
$$
and
$$
f_{n+1}\Bigl(f_{n+1}\bigl(C,\sim\mu_n(b_n)\bigr),\mu_n(b_n)\Bigr)=
f_{n+1}\Bigl(f_{n+1}\bigl(C,\sim\mu_n(b_n)\bigr),\sim\mu_n(b_n)\Bigr)=
f_{n+1}\bigl(C,\sim\mu_n(b_n)\bigr)
\;.
$$
\begin{description}
\item[Proof.]
The result is derived for $f_{n+1}\bigl(C,\mu_n(b_n)\bigr)$\,; it is quite similar for $f_{n+1}\bigl(C,\sim\mu_n(b_n)\bigr)$\,.
\\[3pt]
By corollary 2 of lemma, $f_{n+1}\bigl(C,\mu_n(b_n)\bigr)=\bigl(C\cap\mu_n(b_n)\bigr)\cup\bigl(T(C)\cap\sim\mu_n(b_n)\bigr)$\,.\\
Since $T\bigl(f_{n+1}\bigl(C,\mu_n(b_n)\bigr)\bigr)=f_{n+1}\bigl(C,\mu_n(b_n)\bigr)$ by definition, the proof is done by corollary 2.
\item[$\Box\Box\Box$]\rien
\end{description}
\emph{Corollary.}
As a direct consequence, $f_{n+1}\bigl(f_{n+1}(B,A),A\bigr)=f_{n+1}\bigl(f_{n+1}(B,A),\sim A\bigr)=f_{n+1}(B,A)$\,, whenever $f_{n+1}(B,A)$ exists.
\subsection{Proof of $\tilde\beta5w$}
Assume $f_{n+1}(B,A)$ and $f_{n+1}(B,\sim A)$ exist and $f_{n+1}(B,A)=B$\,.\\
By corollary of lemma 2, $f_{n+1}(B,\sim A)=f_{n+1}\bigl(f_{n+1}(B,A),\sim A\big)=f_{n+1}(B,A)=B$\,.
\section{Proof: completeness for the conditional operator}
\label{Appendix:ProofOfAlmostCompletude}
To be proved:\\[5pt]
\emph{Let $\phi\in\mathcal{L}$ be constructed without $\Box$ or $\Diamond$\,.
Then $\vdash\phi$ in DmBL$_\ast$ if and only if $H_{\mathcal{B}}(\phi)=\Omega_\infty$\,.}
\\[5pt]
From now on, let $\mathcal{L}_b=\bigl\{\phi\in\mathcal{L}\,/\, \phi\mbox{ is constructed without }\Box\mbox{ or }\Diamond\bigr\}$\,.\\
In fact, it will be proved:
\begin{equation}
\label{the:Property}
\mbox{$(H_{\mathcal{B}})_\equiv$ is a Boolean isomorphism between $(\mathcal{L}_b)_\equiv$ and $M_\infty$\,,}
\end{equation}
where $(\mathcal{L}_b)_\equiv$ is the set of equivalence classes of $\mathcal{L}_b$ and $(H_{\mathcal{B}})_\equiv$ is inferred from $H_{\mathcal{B}}$\,.
\\[5pt]
The proof is based on a recursive construction of $\mathcal{L}_b$ similar to the definition of $M_\infty$\,.
\paragraph{Construction.}
Assume the sequence $(\Omega_n,M_n,h_n,f_n,\Lambda_n,\mu_n)_{n\in\Nset}$ being constructed.\\
The sequence $(L_n)_{n\in\Nset}$ is defined by:
\begin{itemize}
\item $L_0=\mathcal{L}_C$\,,
\item $L_{n+1}\subset\mathcal{L}_b$ is the set generated by $L_n$, the classical operators, the conditionals $(\cdot|\phi)$ and $(\cdot|\neg\phi)$ where $\phi\in L_n$ and $H_{\mathcal{B}}(\phi)=b_n$\,.
\end{itemize}
A set $\Sigma_n\subset (L_n)_\equiv$ is called a generating partition of $L_n$, if it verifies:
$$
\forall\phi\in (L_n)_\equiv\,,\;\exists S\subset\Sigma_n\,,\;\bigvee_{\sigma\in S}\sigma=\phi
\quad\mbox{and}\quad
\sigma\wedge\sigma'=\bot\mbox{ for any }\sigma,\sigma'\in\Sigma_n\mbox{ such that }\sigma\ne\sigma'\,.
$$
The following property is proved recursively in the next paragraphs:
\begin{equation}
\label{the:Equation}
\mbox{There is a generating partition }\Sigma_n\mbox{ of }L_n\mbox{ such that }\mathrm{card}(\Sigma_n)\le\mathrm{card}(\Omega_n)\;.
\end{equation}
Since $(H_{\mathcal{B}})_\equiv$ is by construction an onto morphism from $(L_n)_\equiv$ to $M_{n:\infty}$\,,
(\ref{the:Equation}) implies that $(H_{\mathcal{B}})_\equiv$ is a Boolean isomorphism between $(L_n)_\equiv$ and $M_{n:\infty}$\,.\\
The cyclic definition of $\Lambda_n$ then implies $\mathcal{L}_b=\cup_{n\in\Nset}(L_n)_\equiv$ and (\ref{the:Property}) is deduced.
\paragraph{Proof of (\ref{the:Equation}) for $n=0$.}
It is obvious, since $(M_{0:\infty},H_{\mathcal{B}})$ is a complete model for $\mathcal{L}_C$\,.
\paragraph{True for $n$ implies true for $n+1$.}\rien\\
The recursion hypothesis implies that $(H_{\mathcal{B}})_\equiv$ is an isomorphism between $(L_n)_\equiv$ and $M_{n:\infty}$\,.\\
Define then $\beta_n\in (L_n)_\equiv$ such that $(H_{\mathcal{B}})_\equiv(\beta_n)=b_n$\,.\\
It is known that $\bigl((\cdot|\beta_n)\big|\neg\beta_n\bigr)=(\cdot|\beta_n)$ and
$\bigl((\cdot|\neg\beta_n)\big|\beta_n\bigr)=(\cdot|\neg\beta_n)$\,.\\
Then, since sub-universes are classical, $\Sigma_{n+1}=\bigl\{\sigma\wedge(\sigma'|\beta_n)\wedge(\sigma''|\neg\beta_n)\;/\;\sigma,\sigma',\sigma''\in\Sigma_n\bigr\}\setminus\{\bot\}$\,.\\[5pt]
Now, denote $B_n
=\bigl\{\sigma\in\Sigma_n\,\big/\,\sigma\wedge \beta_n=\sigma\bigr\}$
and
$\overline{B}_n
=\bigl\{\sigma\in\Sigma_n\,\big/\,\sigma\wedge\neg\beta_n=\sigma\bigr\}$\,.\\
It comes that $(\sigma'|\beta_n)=(\sigma''|\neg\beta_n)=\bot$ for $\sigma'\in\overline{B}_n$ and $\sigma''\in B_n$\,.\\
Moreover $\sigma\wedge(\sigma'|\beta_n)\wedge(\sigma''|\neg\beta_n)=\bot$ for $\sigma\not\in \{\sigma',\sigma''\}$\,;
on the other hand,
$\sigma\wedge(\sigma|\beta_n)=\sigma$ for $\sigma\in B_n$\,, and
$\sigma\wedge(\sigma|\neg\beta_n)=\sigma$ for $\sigma\in \overline{B}_n$\,.
\\[5pt]
Then, the two construction cases of $(\Omega_n,M_n,h_n,f_n,\Lambda_n)_{n\in\Nset}$ are considered:
\subparagraph{Case 1.}
Then, $\Sigma_{n+1}=
\bigcup_{\sigma\in B_n}
\bigcup_{\sigma'\in \overline{B}_n}\bigl\{\sigma\wedge(\sigma'|\neg\beta_n),\sigma'\wedge(\sigma|\beta_n)\bigr\}
$\,, owing to above discussion.
\\
As a consequence, $
\mathrm{card}(\Sigma_{n+1})\le2\mathrm{card}(B_n)\mathrm{card}(\overline{B}_n)=2\mathrm{card}(b_n)\mathrm{card}(\sim b_n)=\mathrm{card}(\Omega_{n+1})\,.$
\subparagraph{Case 0.}
In this case, $\beta_n=\beta_\nu$\,.\\
Define $C_\nu=\{\sigma\in\Sigma_{\nu+1}/\sigma\wedge\beta_\nu=\sigma\}$
and $\overline{C}_\nu=\{\sigma\in\Sigma_{\nu+1}/\sigma\wedge\neg\beta_\nu=\sigma\}$\,.\\
Define also $D[\phi]=\{\sigma\in\Sigma_n/\sigma\wedge\phi=\sigma\}$ for any $\phi\in(L_{\nu+1})_\equiv$\,.\\[5pt]
From previously, it is know that $\Sigma_{n+1}$ contains elements of the form $\sigma\wedge(\sigma'|\neg\beta_n)$ or $\sigma'\wedge(\sigma|\beta_n)$ with $(\sigma,\sigma')\in B_n\times \overline{B}_n$\,;
but the construction at step $\nu+1$ implies additional constraints, to be specified.
\\[3pt]
Let consider especially the case $\sigma\wedge(\sigma'|\neg\beta_n)$; case $\sigma'\wedge(\sigma|\beta_n)$ is quite similar.\\
Notice that there is $(\tau,\tau')\in C_\nu\times\overline{C}_\nu$
such that $\sigma\in D[\tau]\cap D[(\tau'|\neg\beta_\nu)]$\,,
and $(\theta,\theta')\in \overline{C}_\nu\times C_\nu$
such that $\sigma'\in D[\theta]\cap D[(\theta'|\beta_\nu)]$\,.\\
Now, $\tau\wedge(\tau'|\neg\beta_\nu)\wedge\bigl(\theta\wedge(\theta'|\beta_\nu)\big|\neg\beta_\nu\bigr)=
\bigl(\tau\wedge(\theta'|\beta_\nu)\bigr)\wedge
(\tau'\wedge\theta|\neg\beta_\nu)=\bot$
unless $\tau=\theta'$ and $\tau'=\theta$\,.\\
As a consequence, it is deduced:
$$\sigma\wedge(\sigma'|\neg\beta_n)\ne\bot\quad\mbox{implies}\quad\exists(\tau,\theta)\in C_\nu\times\overline{C}_\nu\,,\;(\sigma,\sigma')\in \bigl(D[\tau]\cap D[(\theta|\neg\beta_\nu)]\bigr)\times\bigl(D[\theta]\cap D[(\tau|\beta_\nu)]\bigr)\;.$$
Similarly, it is deduced:
$$\sigma'\wedge(\sigma|\beta_n)\ne\bot\quad\mbox{implies}\quad\exists(\theta,\tau)\in \overline{C}_\nu\times C_\nu\,,\;(\sigma',\sigma)\in \bigl(D[\theta]\cap D[(\tau|\beta_\nu)]\bigr)\times\bigl(D[\tau]\cap D[(\theta|\neg\beta_\nu)]\bigr)\;.$$
At last $\mathrm{card}(\Sigma_{n+1})\le\sum_{(\tau,\theta)\in C_\nu\times\overline{C}_\nu} 2\,\mathrm{card}\bigl(D[\tau]\cap D[(\theta|\neg\beta_\nu)]\bigr) \mathrm{card}\bigl(D[\theta]\cap D[(\tau|\beta_\nu)]\bigr)$\\
\rien\hspace{175pt}$=\sum_{i\in I_n}2\,\mathrm{card}\bigl(\Pi_n(i)\bigr) \mathrm{card}\bigl(\Gamma_n(i)\bigr)=\mathrm{card}(\Omega_{n+1})\;.$
\section{Probability extension}
\label{Appendix:Probabilition}
To be proved:\\[5pt]
Let $\pi$ be a probability defined over $C$\,, such that $\pi(\phi)>0$ for any $\phi\not\equiv_C\bot$.
Then, there is a (multiplicative) probability $\overline{\pi}$ defined over DmBL$_\ast$ such that  $\forall\phi\in\mathcal{L}_C\,,\;\overline{\pi}(\phi)=\pi(\phi)$\,.
\\[5pt]
The construction of $\overline{\pi}$ is a recursion based on the definition of $(\Omega_n,M_n,h_n,f_n,\Lambda_n,\mu_n)_{n\in\Nset}$\,.
\subsection{Construction}
The probabilities $P_n|_{n\in\Nset}$ are defined over $M_{n:\infty}$ by:
$$
P_n(A_\infty)=\sum_{\omega\in A}P_n\bigl(\omega_\infty\bigr)\quad\mbox{for any }A\in M_n\,,
$$
and:
\subparagraph{Initialization.}\rien\\
For $\omega=\bigl(\delta_\theta|_{\theta\in\Theta}\bigr)\in\Omega_0$ and $\tau\in\Theta$, define $\tau_\omega=\tau$ if $\delta_\tau=1$ and $\tau_\omega=\neg\tau$ if $\delta_\tau=0$\,.\\[3pt]
Then set $P_0\bigl(\omega_\infty\bigr)=\pi\left(\bigwedge_{\theta\in\Theta}\theta_\omega\right)$ for any $\omega\in\Omega_0$\,.
\subparagraph{From n to n+1.}
For any $(\omega,\omega')\in\Pi_n(i)\times\Gamma_n(i)$\,, set:
$$
P_{n+1}\bigl((\omega,\omega')_\infty\bigr)=\frac{P_n(\omega_\infty)P_n(\omega'_\infty)}{P_n\bigl(\Gamma_n(i)_\infty\bigr)}
\quad\mbox{and}\quad
P_{n+1}\bigl((\omega',\omega)_\infty\bigr)=\frac{P_n(\omega_\infty)P_n(\omega'_\infty)}{P_n\bigl(\Pi_n(i)_\infty\bigr)}\;.
$$
An example of construction is given in appendix~\ref{BayesMod:Comp&Ex}.
\subparagraph{Notation.}
For $m\le n$ and $A\in M_m$, the probability $P_n(A_\infty)$ is denoted $P_n(A)$ for simplicity.
\subsection{Properties.}
\subsubsection{Proposition 1} $P_n\subset P_{n+1}$\,, \emph{i.e.} $P_{n+1}(A)=P_{n}(A)$ for any $A\in M_n$\,.
\begin{description}
\item[Proof.]
For $A\in M_n$\,,
$\mu_n(A)=\bigcup_{i\in I_n}\biggl(\Bigl(\bigl(A\cap\Pi_n(i)\bigr)\times\Gamma_n(i)\Bigr)\cup
\Bigl(\bigl(A\cap\Gamma_n(i)\bigr)\times\Pi_n(i)\Bigr)\biggr)$ and then:
$$\rien\hspace{-20pt}\begin{array}{@{}l@{}}\displaystyle
P_{n+1}(A)=\sum_{i\in I_n}\left(
\sum_{\omega\in A\cap\Pi_n(i)}\left(\sum_{\omega'\in\Gamma_n(i)}
\frac{P_n(\omega)P_n(\omega')}{P_n\bigl(\Gamma_n(i)\bigr)}\right)
+
\sum_{\omega\in A\cap\Gamma_n(i)}\left(\sum_{\omega'\in\Pi_n(i)}
\frac{P_n(\omega)P_n(\omega')}{P_n\bigl(\Pi_n(i)\bigr)}\right)
\right)
\vspace{4pt}\\\displaystyle
\rien\hspace{100pt}
=\sum_{i\in I_n}\left(
\sum_{\omega\in A\cap\Pi_n(i)}P_n(\omega)+
\sum_{\omega\in A\cap\Gamma_n(i)}P_n(\omega)
\right)=
\sum_{\omega\in A}P_n(\omega)=P_n(A)\;.
\end{array}$$
\item[$\Box\Box\Box$]\rien
\end{description}
\emph{Corollary.}
$P_n(\Omega_\infty)=1$\,.
\\[3pt]
Derived from $P_0(\Omega_\infty)=\pi(\top)=1$ which is obvious.
\\[4pt]
\emph{Corollary of the corollary.}
$P_n$ is indeed a probability in the classical meaning.
\\[3pt]
Additivity, coherence are obtained by construction.
Finiteness comes from the corollary.
\subsubsection{Proposition 2}
\label{prop:3}
\begin{enumerate}
\item $\displaystyle P_n\bigl(\Pi_n(i)\bigr)+P_n\bigl(\Gamma_n(i)\bigr)=\frac{P_n\bigl(\Pi_n(i)\bigr)}{P_n(b_n)}=\frac{P_n\bigl(\Gamma_n(i)\bigr)}{P_n(\sim b_n)}$\,,
for any $i\in I_n$\,, 
\item
$P_{n+1}\bigl(\mu_n(b_n)\cap A\bigr)=P_{n+1}(b_n)P_{n+1}\Bigl(f_{n+1}\bigl(A,\mu_n(b_n)\bigr)\Bigr)$\,, for any $A\in M_{n+1}$\,,
\item
$P_{n+1}\bigl(\sim \mu_n(b_n)\cap A\bigr)=P_{n+1}(\sim b_n)P_{n+1}\Bigl(f_{n+1}\bigl(A,\sim \mu_n(b_n)\bigr)\Bigr)$\,, for any $A\in M_{n+1}$\,.
\end{enumerate}
These propositions are proved recursively.
\begin{description}
\item[Proof of 1.]
Obvious in case 1; the difficulty arises for case 0.
\\[3pt]
Assume now case 0, and let $(\omega,\omega')\in I_n$\,, \emph{i.e.} $\omega\in\mu_\nu(b_\nu)$ and $\omega'\in\sim\mu_\nu(b_\nu)$.\\
Then $\frac{P_n\Bigl(f_{\nu+1}\bigl(\omega',\sim \mu_\nu(b_\nu)\bigr)\cap\omega\Bigr)}{P_n(b_\nu)}=P_n\biggl(f_{\nu+1}\Bigl(f_{\nu+1}\bigl(\omega',\sim \mu_\nu(b_\nu)\bigr)\cap\omega,\mu_\nu(b_\nu)\Bigr)\biggr)$, by the recursion hypothesis over 2,
and finally $\frac{P_n\bigl(\Pi_n(i)\bigr)}{P_n(b_n)}=P_n\Bigl(f_{\nu+1}\bigl(\omega',\sim \mu_\nu(b_\nu)\bigr)\cap f_{\nu+1}\bigl(\omega,\mu_\nu(b_\nu)\bigr)\Bigr)$\,.\\
Similarly, it is derived $\frac{P_n\bigl(\Gamma_n(i)\bigr)}{P_n(\sim b_n)}=P_n\Bigl(f_{\nu+1}\bigl(\omega',\sim \mu_\nu(b_\nu)\bigr)\cap f_{\nu+1}\bigl(\omega,\mu_\nu(b_\nu)\bigr)\Bigr)$\,.\\
Then $\frac{P_n\bigl(\Pi_n(i)\bigr)}{P_n(b_n)}=\frac{P_n\bigl(\Gamma_n(i)\bigr)}{P_n(\sim b_n)}$ and the result is deduced from $P_n(b_n)+P_n(\sim b_n)=1$.
\item[Proof of 2.]
Since $P_{n+1}\bigl(T(\omega)\bigr)=P_{n+1}(\omega)\frac{P_n\bigl(\Gamma_n(i)\bigr)}{P_n\bigl(\Pi_n(i)\bigr)}$ for $\displaystyle\omega\in \bigcup_{i\in I_n}\bigl(\Pi_n(i)\times\Gamma_n(i)\bigr)$\,,
it comes:
$$\begin{array}{@{}l@{}}
\displaystyle P_{n+1}\Bigl(f_{n+1}\bigl(A,\mu_n(b_n)\bigr)\Bigr)=\sum_{i\in I_n}\biggl(\sum_{\omega\in A\cap(\Pi_n(i)\times\Gamma_n(i))}P_{n+1}(\omega)\frac{P_n\bigl(\Pi_n(i)\bigr)+P_n\bigl(\Gamma_n(i)\bigr)}{P_n\bigl(\Pi_n(i)\bigr)}\biggr)
\vspace{4pt}\\\displaystyle
\rien\hspace{50pt}=\sum_{i\in I_n}\frac{P_{n+1}\Bigl(A\cap\bigl(\Pi_n(i)\times\Gamma_n(i)\bigr)\Bigr)}{P_n(b_n)}
=\frac{P_{n+1}\bigl(\mu_n(b_n)\cap A\bigr)}{P_n(b_n)}\;.
\end{array}$$
\item[Proof of 3.] Similar to 2.
\item[$\Box\Box\Box$]\rien
\end{description}
\subsubsection{Conclusion.}
\label{AppC:Conclude}
Define $P_\infty=\bigcup_{n\in\Nset}P_n$\,, that is $\forall n\in\Nset,\,\forall A\in M_{n:\infty}\,,\;P_\infty(A)=P_n(A)$\,.
\\
By inheritance from $P_n$, $P_\infty$ is a probability over $M_{\infty}$, which verifies the property:
$$
P_{\infty}(A\cap B)=P_{\infty}(A)P_{\infty}\bigl(f_{\infty}(B,A)\bigr)\,,\quad\mbox{for any }A,B\in M_{\infty}\;.
$$
Define $\overline{\pi}(\phi)=P_\infty\bigl(H_{\mathcal{B}}(\phi)\bigr)$\,;
the additivity, coherence and finiteness of $P_\infty$ are inherited by $\overline{\pi}$\,.\\
It is also deduced:
$$
\overline{\pi}(\phi\wedge\psi)=P_\infty\bigl(H_{\mathcal{B}}(\phi)\cap H_{\mathcal{B}}(\psi)\bigr)= P_\infty\bigl(H_{\mathcal{B}}(\phi)\bigr)P_\infty\Bigl(f_{\infty}\bigl(H_{\mathcal{B}}(\psi),H_{\mathcal{B}}(\phi)\bigr)\Bigr)=
\overline{\pi}(\phi)\overline{\pi}\bigl((\psi|\phi)\bigr)\,.
$$
Then, $\vdash\phi\times\psi$ implies $(\psi|\phi)\equiv\psi$ and finally $\overline{\pi}(\phi\wedge\psi)=\overline{\pi}(\phi)\overline{\pi}(\psi)\,.$
\\
At last, $\overline{\pi}$ verifies the multiplicativity.
\\[5pt]
\underline{$\overline{\pi}$ is a (multiplicative) probability over DmBL$_\ast$\,.}
\paragraph{Rational structure of $\overline{\pi}$\,.}
Let $\Sigma=\left\{\left.\bigwedge_{\theta\in\Theta}\epsilon_\theta\;\right/\;\epsilon\in\prod_{\theta\in\Theta}\{\theta,\neg\theta\}\right\}$\,.
For any $\phi\in\mathcal{L}$\,, there is a rational function $R_\phi: \Rset^{(2^\Theta)}\rightarrow\Rset$ such that $\overline{\pi}(\phi)=R_\phi\bigl(\pi(\sigma)|_{\sigma\in\Sigma}\bigr)$\,.
\\[5pt]
The proof is obvious from the construction.
\section{Conditional model:
first steps of construction}
\label{BayesMod:Comp&Ex}
%
%
In this paragraph, the objects $\Omega_{k}, f_{k},\Lambda_{k},\mu_{k-1}|_{k=0,1}$\,, \emph{i.e.} one iteration,  are explicitly constructed, as well as the associated probability extensions $P_0, P_1$ (\emph{c.f.} appendix~\ref{Appendix:Probabilition}).
It is assumed that $\Omega_0=\{a,b,c\}$\,.
This hypothesis cannot hold actually, but the case is sufficiently simple to be handled, and sufficiently complex to be illustrative.
Only the case 1 of the construction is considered.
Case 0 is intractable in a true example. \emph{For simplicity, $\mu_k(A)$ and $A$ are identified.}
\subparagraph{k=0.}
By definition, the list $\Lambda_0$ contains the elements of $\mathcal{P}(\Omega_0)\setminus\{\emptyset,\Omega_0\}$\,.
In this example, it is chosen $\Lambda_0=\{a,b\},c,\{b,c\},a,\{c,a\},b$ and $P_0(a)=0.2$\,, $P_0(b)=0.3$\,, $P_0(c)=0.5$\,.
\subparagraph{k=1.}
Case 1 holds with $\Pi_0(b_0)=b_0=\{a,b\}$ and $\Gamma_0(b_0)=\{c\}$\,.\\
It comes $\Omega_1=\bigl\{(a,c),(b,c),(c,a),(c,b)\bigr\}$\,,
$\mu_0(a)=(a,c)$\,, $\mu_0(b)=(b,c)$\,, $\mu_0(c)=\{(c,a),(c,b)\}$\,, 
\\[3pt]
$P_1(a)=P_1((a,c))=\frac{P_0(a)P_0(c)}{P_0(\Gamma_0(b_0))}=\frac{0.2\times0.5}{0.5}=0.2$\,, $P_1(b)=P_1((b,c))=\frac{0.3\times0.5}{0.5}=0.3$\,,
\\[3pt]
$P_1((c,a))=\frac{P_0(a)P_0(c)}{P_0(\Pi_0(b_0))}=\frac{0.2\times0.5}{0.2+0.3}=0.2$\,, $P_1((c,b))=\frac{0.3\times0.5}{0.2+0.3}=0.3$\,,
\\[3pt]
$f_1\bigl(a,\{a,b\}\bigr)=(\mathrm{id}\cup T)\bigl((a,c)\bigr)=\bigl\{(a,c),(c,a)\bigr\}$\,,
$f_1\bigl(b,\{a,b\}\bigr)=(\mathrm{id}\cup T)\bigl((b,c)\bigr)=\bigl\{(b,c),(c,b)\bigr\}$\,,\\
$f_1\bigl((c,a),c\bigr)=(\mathrm{id}\cup T)\bigl((c,a)\bigr)=\bigl\{(a,c),(c,a)\bigr\}$\,,
$f_1\bigl((c,b),c\bigr)=(\mathrm{id}\cup T)\bigl((c,b)\bigr)=\bigl\{(b,c),(c,b)\bigr\}$\,,
$f_1(c,c)=(\mathrm{id}\cup T)\bigl(\{(c,a),(c,b)\}\bigr)=\Omega_1$\,,
(other cases are obvious)
\\[3pt]
and $\Lambda_1=\{b,c\},a,\{c,a\},b\ +\ \mathcal{P}(\Omega_1)\setminus\mathcal{P}(\Omega_0)\ +\
\{a,b\},c$\,.
\\[5pt]
The relation $P_1\bigl(f_1(B,A)\bigr)P_1(A)=P_1(A\cap B)$ is easily verified for $B\subset\Omega_1$ and $A\in\{b_0,\sim b_0\}\,$.
\end{document}